\documentclass[12pt]{article}
\pdfoutput=1

\usepackage{amsmath}
\usepackage{amssymb,amsfonts}
\usepackage{bm}
\usepackage[mathcal]{euscript}
\usepackage[letterpaper]{geometry}
\usepackage{color}
\usepackage{listings}
\usepackage{natbib}
\usepackage{graphicx}
\usepackage{subfigure}
\usepackage{hyperref}

\usepackage{makeidx}
\makeindex

\graphicspath{{figs/}{figs_lo/}}

\definecolor{beige}{RGB}{245,245,220}

\newcommand{\jlt}[1]{}

\newcommand{\braidlab}{\texttt{braidlab}}
\newcommand{\braid}{\texttt{braid}}
\newcommand{\annbraid}{\texttt{annbraid}}
\newcommand{\loopc}{\texttt{loop}}

\newcommand{\mathnotation}[2]{\newcommand{#1}{\ensuremath{#2}}}

\DeclareMathOperator{\FTBE}{FTBE}

\renewcommand{\l}{\left}                        
\renewcommand{\r}{\right}                       
\mathnotation{\ee}{\mathrm{e}}                  
\mathnotation{\nn}{n}                           
\mathnotation{\ac}{a}                           
\mathnotation{\bc}{b}                           
\mathnotation{\cc}{c}                           
\mathnotation{\dc}{d}                           
\mathnotation{\fc}{f}                           
\mathnotation{\acnew}{\ac'}
\mathnotation{\bcnew}{\bc'}
\mathnotation{\abv}{\bm{u}}                     
\mathnotation{\ip}{i}                           
\newcommand{\pos}[1]{#1^+}
\renewcommand{\neg}[1]{#1^-}
\mathnotation{\ldef}{\mathrel{\raisebox{.069ex}{:}\!\!=}}
\mathnotation{\rdef}{\mathrel{=\!\!\raisebox{.069ex}{:}}}

\begin{document}

\lstset{language=Matlab}
\lstset{breaklines=true}
\lstset{backgroundcolor=\color{beige}}

\lstset{
basicstyle=\small\ttfamily,
keywordstyle=\small\ttfamily,
identifierstyle=,
commentstyle=\small\rmfamily\itshape,
stringstyle=\small\ttfamily,
showstringspaces=false}

\title{\href{{http://github.com/jeanluct/braidlab}}{\braidlab}:
  a software package \\ for braids and loops}
\author{%
\href{http://www.math.wisc.edu/~jeanluc}{Jean-Luc Thiffeault}
  and
\href{http://mbudisic.wordpress.com}{Marko Budi\v{s}i\'{c}}
}
\date{{\small Department of Mathematics\\ University of
    Wisconsin -- Madison} \\[10pt]
  release 3.2.4
} 

\maketitle

\begin{abstract}
  \braidlab\ is a Matlab package for analyzing data using braids.  It was
  designed to be fast, so it can be used on relatively large problems.  It
  uses the object-oriented features of Matlab to provide a class for braids on
  punctured disks and a class for equivalence classes of simple closed loops.
  The growth of loops under iterated action by braids is used to compute the
  topological entropy of braids, as well as for determining the equality of
  braids.  This guide is a survey of the main capabilities of \braidlab, with
  many examples; the help messages of the various commands provide more
  details.  Some of the examples contain novel observations, such as the
  existence of cycles of the linear effective action for arbitrary braids.
\end{abstract}

\tableofcontents

\section{What are braids?}

\subsection{A brief introduction to the braid group}
\label{sec:braidgroup}

\index{braid!group|(}%
\index{braid!geometric|(}%
Braids are collections of strings anchored between two planes,
\begin{figure}
  \begin{center}
    \subfigure[]{
      \includegraphics[height=.2\textheight]{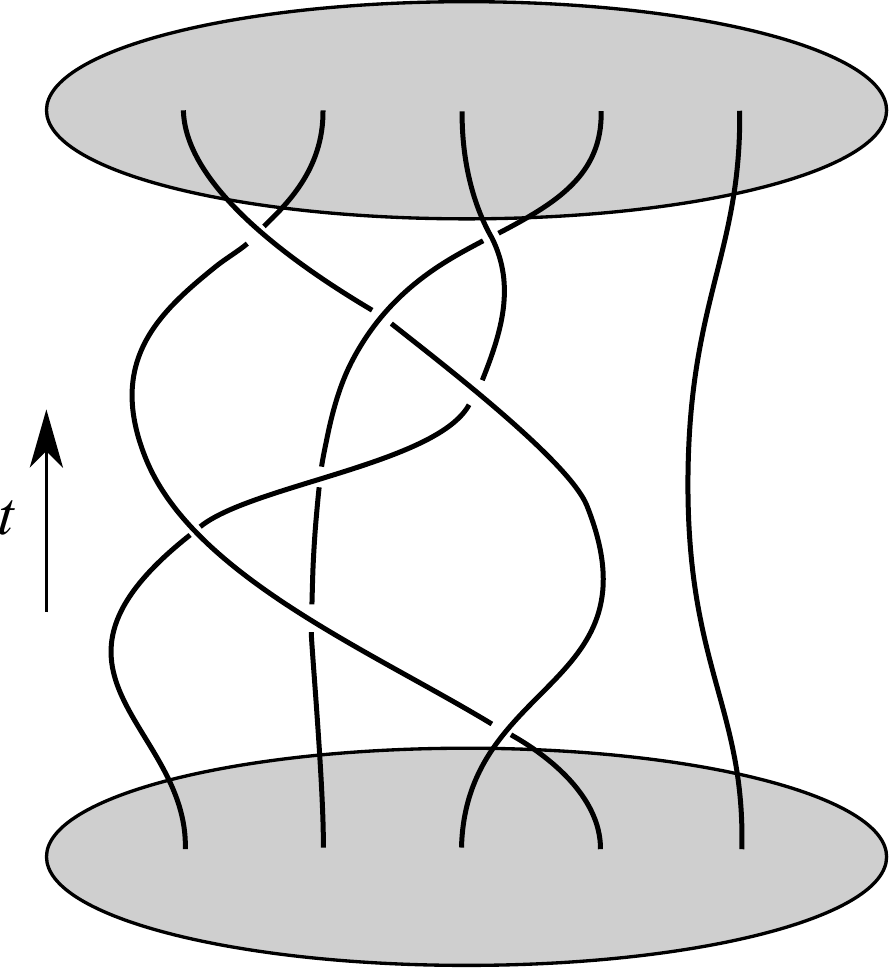}
      \label{fig:geombraid}
    }\hspace{.25\textwidth}
    \subfigure[]{
      \includegraphics[height=.2\textheight]{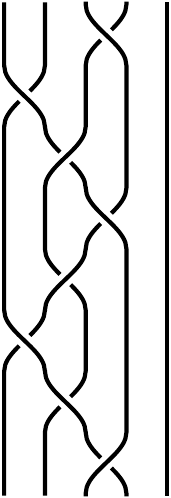}
      \label{fig:sbdiag}
    }
  \end{center}
  \caption{(a) A geometric braid with five strings.  (b) The corresponding
    braid diagram.}
\end{figure}
as in Fig.~\ref{fig:geombraid}.  More precisely this is called a
\emph{geometric braid}.  In this case we can think of the vertical axis as
`time' and the horizontal axes as `space.'  The strings cannot occupy the same
point in space at the same time, and they can't reverse direction (i.e., they
can't go back in time).  We consider two braids to be equal if their strings
can be deformed into each other, with no string crossing another.  The points
on the bottom plane where the strings emanate in Fig.~\ref{fig:geombraid} are
the same as the points on the top plane.  In that sense braids naturally
represent \emph{periodic orbits} of two-dimensional dynamical systems.  When
slicing the braid horizontally at any given time we get a collection of points
in the 2D plane, corresponding to the strings.  We call these points
\emph{punctures} or \emph{particles}. %
\index{punctures}

In fact the set of all braids with a given number of strings, and the same
anchorpoints, form a \emph{group}.  The group multiplication law is simply to
lay one braid after another; it is easy to see that this is associative.  The
identity braid consists of straight strings that are unentangled with each
other.  The inverse of a braid is obtained by reversing time.  (See the book
by \citeauthor{Birman1975} for more details.)  \index{braid!geometric|)}%

\index{braid!algebraic|(}%
A convenient way to represent a braid is to deform it %
as in Fig,~\ref{fig:sbdiag}.  Here the geometric braid is combed such that
only one crossing occurs at a time, and each crossing fits in a time interval
of the same length.  Such a picture is called a \emph{braid diagram}.  The
operation of the~$i$th string (counted from left to right) being exchanged
with the~$(i+1)$th string, such that the left string passes over the right, is
called a %
\index{braid!generator|(}%
\emph{generator} of the braid group, as is denoted~$\sigma_i$.  Any element of
the braid group on~$\nn$ strings can be written as a product of the
generators~$\{\sigma_1,\ldots,\sigma_{\nn-1}\}$ and their inverses.  Thus, the
braid in Fig.~\ref{fig:sbdiag} can be written
\begin{equation}
  \sigma_3\,\sigma_2^{-1}\sigma_1^{-1}\,\sigma_2\sigma_3^{-1}
  \sigma_2\sigma_1^{-1}\sigma_3^{-1}\,,
  \label{eq:sbdiag}
\end{equation}
where we read the generators from left to right, and from bottom to top in
Fig.~\ref{fig:sbdiag}.  A braid written in terms of generators as
in~\eqref{eq:sbdiag} is called an \emph{algebraic braid}.
\index{braid!algebraic|)}%

\begin{figure}
  \begin{center}
    \subfigure[]{
      \raisebox{.09\textwidth}{%
        \includegraphics[width=.35\textwidth]{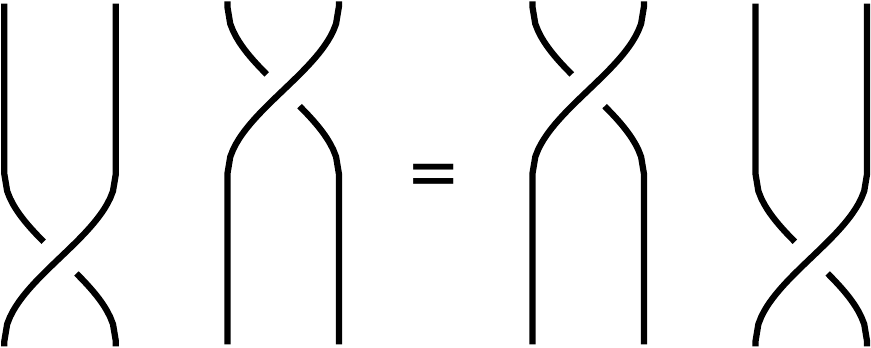}}
      \label{fig:13is31}
    }\hspace{.1\textwidth}
    \subfigure[]{
      \includegraphics[width=.35\textwidth]{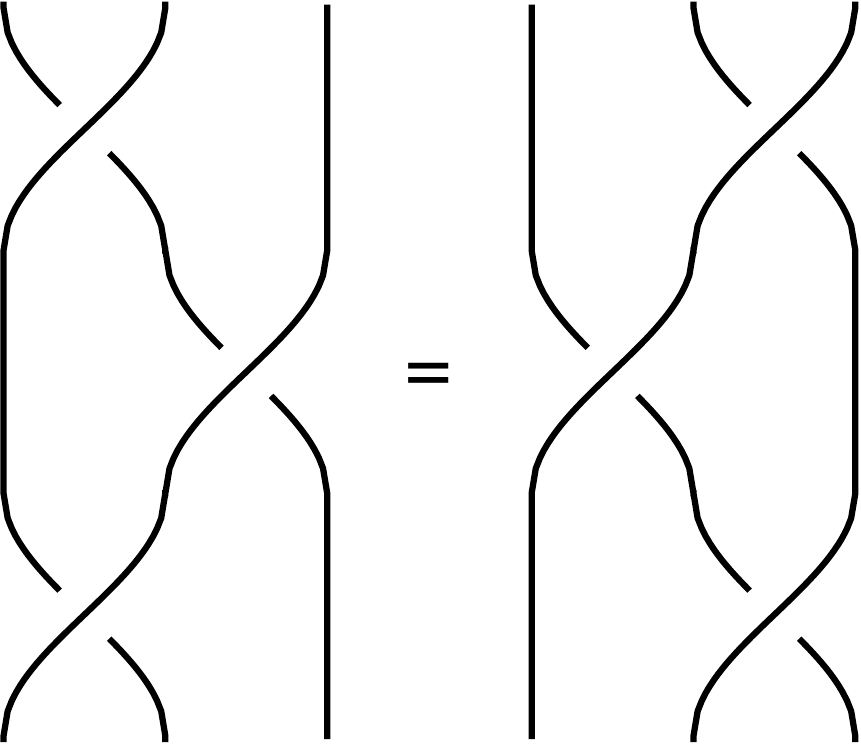}
      \label{fig:121is212}
    }
  \end{center}
  \caption{(a) Generators that don't share a string commute.  (b) The braid
    relation.}
  \label{fig:relations}
\end{figure}
The generators obey some special rules, called \emph{relations}, by virtue of
arising from geometrical braids.  The relations are
\begin{equation}
  \sigma_j\sigma_k=\sigma_k\sigma_j, \quad \lvert j-k\rvert > 1; \qquad
  \sigma_j\sigma_k\sigma_j = \sigma_k\sigma_j\sigma_k, \quad \lvert j-k\rvert=1.
  \label{eq:relations}
\end{equation}
The first type of relation, depicted in Fig.~\ref{fig:13is31}, says that
generators commute if they don't share a string.  The second type, often
called the braid relation, reflects the equality of the two braids shown in
Fig.~\ref{fig:121is212}.  While it's obvious from Fig.~\ref{fig:relations}
that the relations~\eqref{eq:relations} are satisfied, it is far less obvious
that those are the \emph{only} relations that hold, as proved by
\citet{Artin1947}.  Hence, the relations~\eqref{eq:relations} fully
characterize the braid group for~$\nn$ strings, denoted~$B_\nn$.

\index{braid!generator|)}
\index{braid!group|)}%

\subsection{Constructing a braid from orbit data}
\label{sec:braidfromorbitdata}

\index{braid!from data|(}%
\begin{figure}
  \begin{center}
    \subfigure[]{
      \includegraphics[width=.3\textwidth]{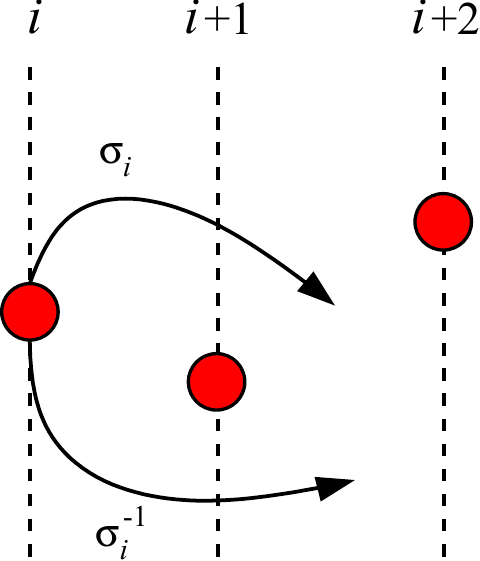}
      \label{fig:crossing}
    }\hspace{.2\textwidth}
    \subfigure[]{
      \includegraphics[width=.3\textwidth]{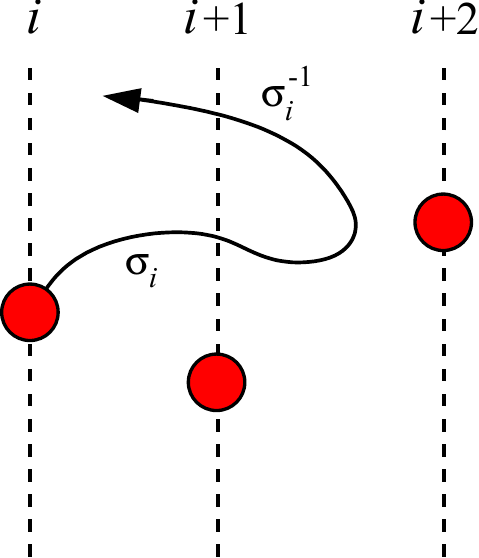}
      \label{fig:crossing2}
    }
  \end{center}
  \caption{(a) The two types of generators that occur when particles exchange
    positions.  (b) A particle exchanging position twice in a row with another
    leads to the generators~$\sigma_i$ followed by~$\sigma_i^{-1}$, which
    cancel each other. [After~\citet{Thiffeault2005}.]}
\end{figure}

So far we have regarded braids as geometrical objects, and showed how to turn
them into algebraic objects by writing them in terms of generators.  But in
practice how do we create braids from particle orbits?  That is, if we have
two-dimensional continuous trajectory data arising from some dynamical system,
how do we turn this data into a braid?

\index{projection line|(}%
The technique to do this was described in \citet{Thiffeault2005}.  We define
an arbitrary line in the 2D plane, called the \emph{projection line}, and look
at the order of the particles projected along that line.  We label each
particle according to its order when projected on this line.  In
Fig.~\ref{fig:crossing}, for example, we see three particles labeled~$i$,
$i+1$, $i+2$, where the projection line is the horizontal.

\index{crossing|(}%
A \emph{crossing} occurs whenever any two particles exchange position (only
adjacent particles can exchange position at a given time, since the
trajectories are continuous).  For particle~$i$ exchanging position with
particle~$i+1$, we assign a generator~$\sigma_i$ to the crossing if the
particle on the left passes above the one on the right, and~$\sigma_i^{-1}$ if
it passes below (Fig.~\ref{fig:crossing}).  As we watch the trajectories
unfold, we construct an algebraic braid as a sequence of generators, with each
generator corresponding to a crossing, and their order determined by when each
crossing occurs.  Note that when two particles exchange position along the
crossing line twice as in Fig.~\ref{fig:crossing2}, without exchanging
positions vertically, then the two crossings yield the generators~$\sigma_i$
followed by~$\sigma_i^{-1}$, which cancel. %
\index{projection line|)}%
\index{crossing|)}%

The outcome of this procedure is not a true geometric braid, since the
particles do not necessarily end up in the same positions as they started.
(If the particles are part of a periodic orbit, then we do obtain a geometric
braid.)  We do however obtain an algebraic braid, and as long as we are
careful to take long enough trajectories the fact that the geometric braid
does not close will not be too consequential.

\index{braid!from data|)}%

\section{A tour of \braidlab}
\label{sec:tour}

You will need access to a recent version of Matlab to use \braidlab.  See
Appendix~\ref{sec:install} for instructions on how to install \braidlab\ on
your machine.

\subsection{The \braid\ class}
\label{sec:braidclass}

\index{braid class@\braid\ class|(}

\subsubsection{Constructor and elementary operations}

\braidlab\ defines a number of classes, most importantly \braid\ and \loopc.
The braid~$\sigma_1\sigma_2^{-1}$ is constructed with %
\index{braid class@\braid\ class!constructor|(}
\begin{lstlisting}[frame=single,framerule=0pt]
>> a = braid([1 -2])   % defaults to 3 strings

a = < 1 -2 >
\end{lstlisting}
which defaults to the minimum required strings,~$3$.  The same braid
on~$4$ strings is constructed with
\begin{lstlisting}[frame=single,framerule=0pt]
> a4 = braid([1 -2],4)   % force 4 strings

a4 = < 1 -2 >
\end{lstlisting}
Two braids can be multiplied: %
\index{braid class@\braid\ class!multiplication (\lstinline{*})}
\begin{lstlisting}[frame=single,framerule=0pt]
>> a = braid([1 -2]); b = braid([1 2]);
>> a*b, b*a

ans = < 1 -2  1  2 >

ans = < 1  2  1 -2 >
\end{lstlisting}
Powers %
\index{braid class@\braid\ class!power (\lstinline{^})}%
can also be taken, including the inverse: %
\index{braid class@\braid\ class!inverse (\lstinline{inv})}%
\begin{lstlisting}[frame=single,framerule=0pt]
>> a^5, inv(a), a*a^-1

ans = < 1 -2  1 -2  1 -2  1 -2  1 -2 >

ans = < 2 -1 >

ans = < 1 -2  2 -1 >
\end{lstlisting}
\index{braid class@\braid\ class!identity braid}%
Note that this last expression is the identity braid, but is not simplified.
The method \lstinline{compact} attempts to simplify the braid: %
\index{braid class@\braid\ class!compact@\lstinline{compact}}
\begin{lstlisting}[frame=single,framerule=0pt]
>> compact(a*a^-1)

ans = < e >
\end{lstlisting}
The method \lstinline{compact} is based on the heuristic algorithm
of~\citet{Bangert2002}, since finding the braid of minimum length in the
standard generators is in general difficult~\citep{Paterson1991}.  Hence,
there is no guarantee that in general \lstinline{compact} will find the
identity braid, even though it does so here.  To really test if a braid is the
identity (trivial braid), use the method \lstinline{istrivial}: %
\index{braid class@\braid\ class!istrivial@\lstinline{istrivial}}
\begin{lstlisting}[frame=single,framerule=0pt]
>> istrivial(a*a^-1)

ans = 1
\end{lstlisting}

The number of strings is %
\index{braid class@\braid\ class!number of strings (\lstinline{n})}
\begin{lstlisting}[frame=single,framerule=0pt]
>> a.n

ans = 3
\end{lstlisting}
Note that
\index{help@\lstinline{help}|(}
\begin{lstlisting}[frame=single,framerule=0pt]
>> help braid
\end{lstlisting}
describes the class \braid.  To get more information on the \braid\
constructor, invoke%
\index{braid class@\braid\ class!constructor}%
\begin{lstlisting}[frame=single,framerule=0pt]
>> help braid.braid
\end{lstlisting}
which refers to the method \braid\ within the class \braid. %
\index{help@\lstinline{help}|)}%
(Use \lstinline{methods(braid)} to list all the methods in the class.)  There
are other ways to construct a \braid, such as using random %
\index{braid class@\braid\ class!constructor!random braid} generators, here a
braid with~$5$ strings and~$10$ random generators:
\begin{lstlisting}[frame=single,framerule=0pt]
>> braid('Random',5,10)

ans = < 1  4 -4  2  4 -1 -2  4  4  4 >
\end{lstlisting}
The constructor can also build some standard braids: %
\index{braid class@\braid\ class!constructor!half-twist}
\index{braid class@\braid\ class!constructor!knots}
\index{knot!braid representative}
\begin{lstlisting}[frame=single,framerule=0pt]
>> braid('HalfTwist',5)

ans = < 4  3  2  1  4  3  2  4  3  4 >

>> braid('8_21')  % braid for 8-crossing knot #21

ans = < 4  3  2  1  4  3  2  4  3  4 >
\end{lstlisting}
In Section~\ref{sec:braidfromdata} we will show how to construct a braid from
a trajectory data set.  \index{braid class@\braid\ class!constructor|)}

\index{braid class@\braid\ class!equality (\lstinline{==})}%
The \braid\ class handles equality of braids:
\begin{lstlisting}[frame=single,framerule=0pt]
>> a = braid([1 -2]); b = braid([1 -2 2 1 2 -1 -2 -1]);
>> a == b

ans = 1
\end{lstlisting}
These are the same braid, even though they appear different from their
generator sequence~\citep{Birman1975}.  Equality is determined efficiently
by %
\index{loop!coordinates}%
\index{Dynnikov coordinates|see{loop coordinates}}%
acting on loop coordinates~\citep{Dynnikov2002}, as described by
\citet{Dehornoy2008}.  See Sections~\ref{sec:loop}--\ref{sec:loopcoords} for
more details.  If for some reason lexicographic (generator-per-generator)
equality of braids is needed, use the method \lstinline{lexeq(b1,b2)}. %
\index{braid class@\braid\ class!lexeq class@\lstinline{lexeq}}

We can extract a subbraid %
\index{braid class@\braid\ class!subbraid class@\lstinline{subbraid}|(} by
choosing specific strings: for example, if we take the~$4$-string
braid~$\sigma_1\sigma_2\sigma_3^{-1}$ and discard the third string, we
obtain~$\sigma_1\sigma_2^{-1}$:
\begin{lstlisting}[frame=single,framerule=0pt]
>> a = braid([1 2 -3]);
>> subbraid(a,[1 2 4])   % subbraid using strings 1,2,4

ans = < 1 -2 >
\end{lstlisting}
\index{braid class@\braid\ class!subbraid class@\lstinline{subbraid}|)}

The opposite of subbraid is the \emph{tensor product}, the larger braid
obtained by laying two braids side-by-side \citep{KasselTuraev}: %
\index{braid class@\braid\ class!tensor@\lstinline{tensor}|(}%
\begin{lstlisting}[frame=single,framerule=0pt]
>> a = braid([1 2 -3]); b = braid([1 -2]);
>> tensor(a,b)

ans = < 1  2 -3  5 -6 >
\end{lstlisting}
Here, the tensor product of a 4-braid and a 3-braid has 7 strings.  The
generators $\sigma_1\sigma_2^{-1}$ of \lstinline{b} became
$\sigma_5\sigma_6^{-1}$ after re-indexing so they appear to the right of
\lstinline{a}.

\subsubsection{Topological entropy and complexity}
\label{sec:entropy}

There are a few methods that exploit the connection between braids and
homeomorphisms \index{homeomorphism} of the punctured disk. %
\index{disk, punctured}%
\index{punctures}
Braids label \emph{isotopy classes} %
\index{homeomorphism!isotopy classes} of homeomorphisms, so we can assign a
topological entropy %
\index{braid class@\braid\ class!entropy@\lstinline{entropy}|(}%
\index{braid!entropy|(}%
\index{topological entropy|see{braid entropy}}%
\index{entropy|see{braid entropy}}%
to a braid:
\begin{lstlisting}[frame=single,framerule=0pt]
>> entropy(braid([1 2 -3]))

ans = 0.8314
\end{lstlisting}
\index{action!of braid on loop}%
The entropy is computed by iterated action on a loop~\citep{Moussafir2006}.
This can fail if the braid is finite-order %
\index{braid!finite-order|(}%
or has very low entropy:
\begin{lstlisting}[frame=single,framerule=0pt]
>> entropy(braid([1 2]))
Warning: Failed to converge to requested tolerance; braid is likely finite-order or has low entropy.  Returning zero entropy.

ans = 0
\end{lstlisting}
To force the entropy to be computed using the Bestvina--Handel train track
algorithm~\citep{Bestvina1995}, %
\index{Bestvina--Handel algorithm|(}%
we add an optional \lstinline{'Method'} parameter:
\begin{lstlisting}[frame=single,framerule=0pt]
>> entropy(braid([1 2]),'Method','train')

ans = 0
\end{lstlisting}
\index{braid!finite-order|)}%
Note that for large braids the Bestvina--Handel algorithm is impractical.

The topological entropy is a measure of braid complexity that relies on %
\index{braid class@\braid\ class!entropy@\lstinline{entropy}|)}%
\index{braid!entropy|)}%
iterating the braid.  It gives the maximum growth rate of a `rubber band'
anchored on the braid, as the rubber band slides up many repeated copies of
the braid.  For finite-order braids, %
\index{braid!finite-order}%
this will converge to zero.  The \emph{geometric complexity} %
\index{braid!complexity|(}%
\index{braid class@\braid\ class!complexity@\lstinline{complexity}|(}%
\index{geometric complexity|see{braid complexity}}%
\index{complexity|see{braid complexity}}%
of a braid~\citep{Dynnikov2007}, is defined in terms of the $\log_2$ of the
number of intersections of a set of curves with the real axis, after one
application of the braid:
\begin{lstlisting}[frame=single,framerule=0pt]
>> complexity(braid([1 -2]))

ans = 2

>> complexity(braid([1 2]))

ans = 1.5850
\end{lstlisting}
See Section~\ref{sec:loop} or `\lstinline{help braid.complexity}' for details
on how the geometric complexity is computed. %
\index{braid!complexity|)}%
\index{braid class@\braid\ class!complexity@\lstinline{complexity}|)}%

\subsubsection{Train track map and transition matrix}

\index{train track map|(}%
\index{transition matrix|(}%
The Bestvina--Handel train track algorithm~\citep{Bestvina1995} can be used to
determine the Thurston--Nielsen type %
\index{braid!Thurston--Nielsen type \lstinline{tntype}}%
of the braid as well as the train track map and its transition
matrix~\citep{Fathi1979,Thurston1988,Casson1988,Boyland1994}: \index{braid
  class@\braid\ class!train@\lstinline{train}|(}%
\index{braid!pseudo-Anosov|(}%
\index{braid!reducible|(}%
\begin{lstlisting}[frame=single,framerule=0pt]
>> train(braid([1 2 -3]))

ans = struct with fields:

       braid: [1x1 braidlab.braid]
      tntype: 'pseudo-Anosov'
     entropy: 0.8314
    transmat: [4x4 double]
       ttmap: {8x1 cell}

>> train(braid([1 2]))

ans = struct with fields:

       braid: [1x1 braidlab.braid]
      tntype: 'finite-order'
     entropy: 0
    transmat: [3x3 double]
       ttmap: {6x1 cell}

>> train(braid([1 2],4))  % reducing curve around 1,2,3

ans = struct with fields:

       braid: [1x1 braidlab.braid]
      tntype: 'reducible'
     entropy: 0
    transmat: [3x3 double]
       ttmap: {7x1 cell}
\end{lstlisting}
\index{braid!pseudo-Anosov|)}%
\index{braid!reducible|)}%
\braidlab\ uses Toby Hall's implementation of the Bestvina--Handel
algorithm~\citep{HallTrain}. %
\index{Bestvina--Handel algorithm|)}%
\index{braid class@\braid\ class!train@\lstinline{train}|)}%

The train track map can be displayed in a human-readable format using the
command \lstinline{ttmap}: %
\index{ttmap@\lstinline{ttmap}|(}%
\begin{lstlisting}[frame=single,framerule=0pt]
>> tt = train(braid([1 2 -3]));
>> ttmap(tt)

 1 -> 4
 2 -> 1
 3 -> 2
 4 -> 3
 a -> D
 b -> d a -3 b -4 B
 c -> B 3 A
 d -> c
\end{lstlisting}
Here peripheral (infinitesimal) edges are denoted by numbers and main edges by
letters.  Inverse main edges are denoted by capital letters.  The display of
infinitesimal edges can be suppressed:
\begin{lstlisting}[frame=single,framerule=0pt]
>> ttmap(tt,'Peripheral',false)

 a -> D
 b -> d a b B
 c -> B A
 d -> c
\end{lstlisting}
The transition matrix associated with the train track map does \emph{not}
contain the peripheral edges, since these do not affect the entropy:
\begin{lstlisting}[frame=single,framerule=0pt]
>> tt.transmat

ans = 0     1     1     0
      0     2     1     0
      0     0     0     1
      1     1     0     0

>> max(abs(eig(ans)))

ans = 2.2966
\end{lstlisting}
\index{transition matrix|)}%
\index{train track map|)}%
\index{ttmap@\lstinline{ttmap}|)}%

\subsubsection{Representation and invariants}

There are a few remaining methods in the braid class, which we describe
briefly.  The reduced Burau matrix
representation~\citep{Burau1936,Birman1975} %
\index{braid!Burau representation|(}%
\index{braid class@\braid\ class!burau@\lstinline{burau}|(}%
of a braid is obtained with the method \lstinline{burau}:
\begin{lstlisting}[frame=single,framerule=0pt]
>> burau(braid([1 -2]),-1)

ans = 1    -1
     -1     2
\end{lstlisting}
where the last argument ($-1$) is the value of the parameter~$t$ in the
Laurent polynomials %
\index{Laurent polynomials|(}%
that appear in the entries of the Burau matrices.  With access to Matlab's
wavelet toolbox, %
\index{Matlab!wavelet toolbox}%
\index{laurpoly@\lstinline{laurpoly}|(}%
we can use actual Laurent polynomials as the entries:
\begin{lstlisting}[frame=single,framerule=0pt]
>> B = burau(braid([1 -2]),laurpoly(1,1))

     | - z^(+1)        z^(+1)     |
     |                            |
 B = |                            |
     |                            |
     |   - 1        + 1 - z^(-1)  |
\end{lstlisting}
but the matrix is now given as a cell array%
\index{Matlab!cell array}%
\index{cell array|see{Matlab cell array}}%
\footnote{%
  \index{Matlab!cell array}%
  A Matlab cell array is similar to a numeric array, except that its entries
  can hold any data, not just numeric.  The entries are indexed as
  \lstinline{a\{1,2\}} rather than \lstinline{a(1,2)}, and matrix operations
  like multiplication are not defined.} %
, %
each entry containing a \lstinline{laurpoly} object:
\begin{lstlisting}[frame=single,framerule=0pt]
>> B{2,2}

ans(z) = + 1 - z^(-1)
\end{lstlisting}
Instead of \lstinline{laurpoly} objects, we can use Matlab's symbolic
toolbox: \index{Matlab!symbolic toolbox}%
\begin{lstlisting}[frame=single,framerule=0pt]
>> B = burau(braid([1 -2]),sym('t'))

B = [ -t,       t]
    [ -1, 1 - 1/t]
\end{lstlisting}
where now \lstinline{B} is a matrix of \lstinline{sym} objects:
\begin{lstlisting}[frame=single,framerule=0pt]
>> B(2,2)

ans = 1 - 1/t
\end{lstlisting}

\index{braid!Lawrence--Krammer representation|(}%
\index{braid class@\braid\ class!lk@\lstinline{lk}|(}%
Another well-known homological representation of braid groups is the
Lawrence--Krammer representation \citep{Lawrence1990,Bigelow2001}.  It is
given in terms of two parameters, usually denoted~$t$ and~$q$:
\begin{lstlisting}[frame=single,framerule=0pt]
>> K = lk(braid([1 -2]),sym('t'),sym('q'))

K = [ -(q-1)^2/q - q*t*(q-1),  -(q-1)/q,  (q^2-q+1)/q^2]
    [                 -q^2*t,         0,              0]
    [           -(q-1)/(q*t),  -1/(q*t),  (q-1)/(q^2*t)]
\end{lstlisting}
In this case there we cannot use \lstinline{laurpoly} entries, since the
representation involves Laurent polynomials in two symbols.  For this reason,
and because its size grows more rapidly with the number of strings (matrices
of dimension $\tfrac12\nn(\nn-1)$), the Lawrence--Krammer representation is
very slow to compute for large braids. %
\index{braid!Lawrence--Krammer representation|)}%
\index{braid class@\braid\ class!lk@\lstinline{lk}|)}%

The reduced Burau matrix of a braid can be used to compute the
\emph{Alexander--Conway polynomial} (or Alexander polynomial for short) %
\index{Alexander--Conway polynomial|(}%
of its closure.  For instance, the trefoil knot is given by the closure of
the %
\index{knot!trefoil}%
\index{knot!Alexander polynomial}%
braid~$\sigma_1^3$ \citep{AlexanderPolynomial}, which gives a Laurent
polynomial
\index{braid class@\braid\ class!alexpoly@\lstinline{alexpoly}|(}%
\begin{lstlisting}[frame=single,framerule=0pt]
>> alexpoly(braid([1 1 1]))  % can also use braid('Trefoil')

ans(z) = + z^(+2) - z^(+1) + 1
\end{lstlisting}
The figure-eight knot is the closure of~$(\sigma_1\sigma_2^{-1})^2$:
\index{knot!figure-eight}
\begin{lstlisting}[frame=single,framerule=0pt]
>> alexpoly(braid([1 -2 1 -2]))  % or braid('Figure-8')

ans(z) = - 1 + 3*z^(-1) - z^(-2)
\end{lstlisting}
This can be `centered' so that it satisfies~$p(z)=\pm p(1/z)$:
\begin{lstlisting}[frame=single,framerule=0pt]
>> alexpoly(braid([1 -2 1 -2]),'Centered')

ans(z) = - z^(+1) + 3 - z^(-1)
\end{lstlisting}
The centered Alexander polynomial is a knot invariant, \index{knot!invariant}
so it can be used to determine when two knots are not the same.  For knots,
the centered polynomial is guaranteed to have integral powers.  For links,
such as the Hopf link consisting of two singly-linked loops, it might not:
\begin{lstlisting}[frame=single,framerule=0pt]
>> alexpoly(braid([1 1]),'Centered')  % the Hopf link

Error using braidlab.braid/alexpoly
Polynomial with fractional powers.  Remove 'Centered' option or use the symbolic toolbox.
\end{lstlisting}
\index{Matlab!symbolic toolbox}%
Fractional powers cannot be represented with a \lstinline{laurpoly} object. %
\index{laurpoly@\lstinline{laurpoly}|)}%
\index{Laurent polynomials|)}%
In that case we can drop the \lstinline{Centered} option, which yields the uncentered
polynomial $1-z$.  Alternatively, we can switch to using a variable from the
symbolic toolbox:
\begin{lstlisting}[frame=single,framerule=0pt]
>> alexpoly(braid([1 1]),sym('x'),'Centered')

ans = 1/x^(1/2) - x^(1/2)
\end{lstlisting}
\index{sym@\lstinline{sym}}%
which can represent fractional powers.  This polynomial
satisfies~$p(x)=-p(1/x)$.

\index{braid!Burau representation|)}%
\index{braid class@\braid\ class!burau@\lstinline{burau}|)}%
\index{Alexander--Conway polynomial|)}%
\index{braid class@\braid\ class!alexpoly@\lstinline{alexpoly}|)}%

The method \lstinline{perm} %
\index{braid class@\braid\ class!perm@\lstinline{perm}|(}%
gives the permutation of strings corresponding to a braid: %
\begin{lstlisting}[frame=single,framerule=0pt]
>> perm(braid([1 2 -3]))

ans = 2  3  4  1
\end{lstlisting}
\index{braid class@\braid\ class!perm@\lstinline{perm}|)}%
If the strings are unpermuted, then the braid is \emph{pure}, %
\index{braid!pure}%
which can also be tested with the method \lstinline{ispure}. %
\index{braid class@\braid\ class!ispure@\lstinline{ispure}}%

Finally, the \emph{writhe} %
\index{braid!writhe}%
\index{braid class@\braid\ class!writhe@\lstinline{writhe}|(}%
of a braid is the sum of the powers of its generators.  The writhe of
$\sigma_1^{+1}\sigma_2^{+1}\sigma_3^{-1}$ is $+1+1-1 = 1$:
\begin{lstlisting}[frame=single,framerule=0pt]
>> writhe(braid([1 2 -3]))

ans = 1
\end{lstlisting}
The writhe is a braid invariant.
\index{braid class@\braid\ class!writhe@\lstinline{writhe}|)}%

\subsubsection{The \annbraid\ subclass}
\label{sec:annbraid}

\index{annbraid@\lstinline{annbraid}|(}
\index{punctures!in annulus|(}
\begin{figure}
  \begin{center}
    \subfigure[]{
      \includegraphics[width=.4\textwidth]{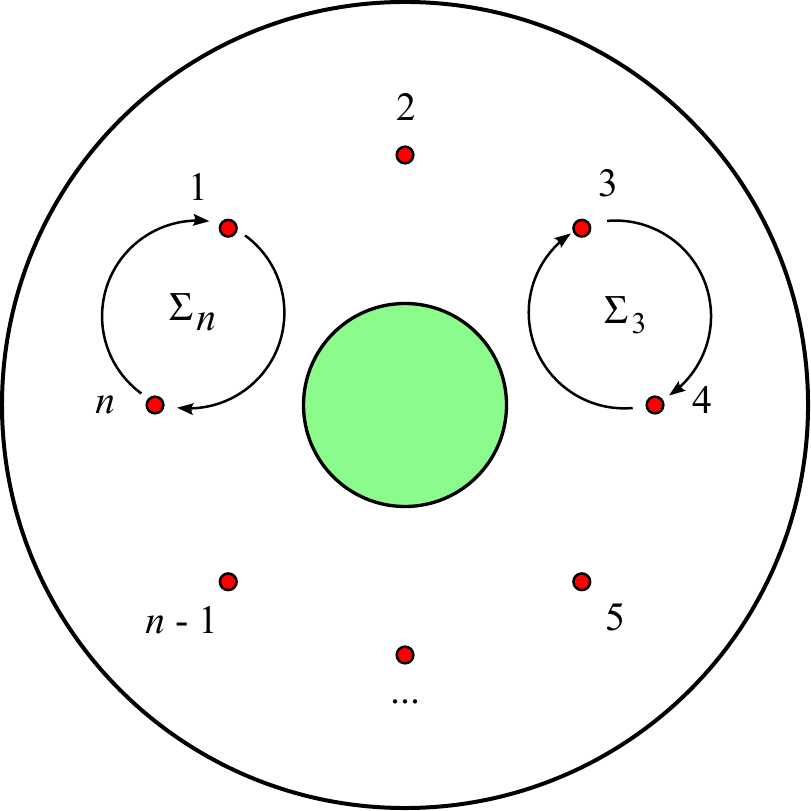}
      \label{fig:annulardomain}
    }\hspace{1em}
    \subfigure[]{
      \raisebox{3.75em}{\includegraphics[width=.5\textwidth]{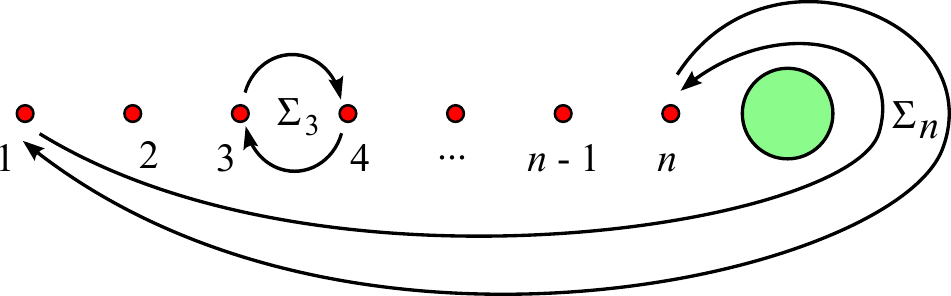}}
      \label{fig:annulardomain2}
    }
  \end{center}
  \caption{(a) Punctures in an annular domain, with two generators.  The
    generator~$\Sigma_\nn$ is unique to the annulus.  (b) The punctures
    rearranged with the center of the annulus as an extra puncture on the
    right, showing how the generator~$\Sigma_\nn$ can be deformed in terms of
    standard generators as in~\eqref{eq:Sigman}.}
\end{figure}

Often it is useful to consider braids in an annular domain, as in
Fig.~\ref{fig:annulardomain}.  It is convenient to rearrange the punctures as
in Fig.~\ref{fig:annulardomain2}, since the overall topology is unchanged.
The center of the annulus becomes an extra puncture, but that extra puncture
is fixed.  For~$\nn$ moving punctures, the braid group on the annulus
has~$\nn$ generators~$\{\Sigma_1,\ldots,\Sigma_\nn\}$, one more than for the
standard braid group, owing to the fixed puncture.  These are related to the
standard generators by~$\Sigma_i=\sigma_i$ for~$1 \le i <\nn$, and
\begin{equation}
  \Sigma_\nn = \sigma_\nn^2 \sigma_{\nn-1} \cdots \sigma_2\,\sigma_1\,
  \sigma_2^{-1} \cdots \sigma_{\nn-1}^{-1} \sigma_\nn^{-2}\,.
  \label{eq:Sigman}
\end{equation}
This last generator effectively exchanges puncture~$\nn$ with puncture~$1$,
exploiting the annular topology.  These kinds of braids were considered in
\citet{Boyland1994,MattFinn2006,MattFinn2011_silver}.

\braidlab\ supports annular braids with the subclass \annbraid, derived from
\braid.  The syntax for creating an annular braid is
\begin{lstlisting}[frame=single,framerule=0pt]
>> b = annbraid([1 2 -3])

b = < 1  2 -3 >*
\end{lstlisting}
The asterisk indicates that this is an annular braid, which has an extra fixed
puncture on the right (called the basepoint).  Hence, the braid has~$4$
punctures, as indicated by
\begin{lstlisting}[frame=single,framerule=0pt]
>> b.n    % total number of punctures, including the fixed one

ans = 4
\end{lstlisting}
but only~$3$ punctures can move, as returned by \lstinline{nann}, the number
of annular punctures:
\begin{lstlisting}[frame=single,framerule=0pt]
>> b.nann    % number of moving punctures

ans = 3
\end{lstlisting}

Many of the methods work described for the \lstinline{braid} class can be
applied to \annbraid{}s.  For instance,
\begin{lstlisting}[frame=single,framerule=0pt]
>> entropy(braid([1 -2])))      % entropy of a normal braid

ans = 0.9624

>> entropy(annbraid([1 -2])))   % entropy of annular braid

ans = 1.7627
\end{lstlisting}
The annular braid has more entropy, since curves grow faster by getting
entangled on the extra puncture \citep{MattFinn2011_silver}.
\begin{figure}
  \begin{center}
    \includegraphics[width=.7\textwidth]{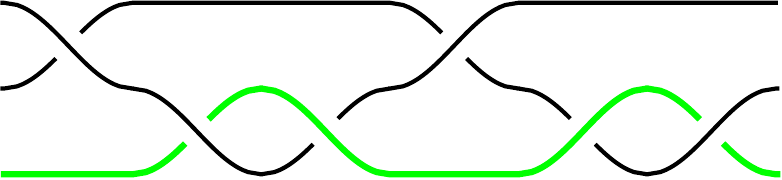}
  \end{center}
  \caption{The output of \lstinline{plot(annbraid([1 -2]))}.  The green strand
    represents the center of the annulus, as in
    Fig.~\ref{fig:annulardomain2}. (See Section~\ref{sec:prop} for how to plot
    braids sideways.)}
  \label{fig:annbraid_s1s-2}
\end{figure}
The annular braid is shown in Fig.~\ref{fig:annbraid_s1s-2}.

\index{punctures!in annulus|)}
\index{annbraid@\lstinline{annbraid}|)}

\subsection{Constructing a braid from data}
\label{sec:braidfromdata}

\subsubsection{An example}
\label{sec:braidfromdataex}

\index{braid class@\braid\ class!constructor!from data|(}%
\index{braid!from data|(}%
One of the main purposes of \braidlab\ is to analyze two-dimensional
trajectory data using braids.  We can assign a braid to trajectory data by
looking for \emph{crossings} %
\index{crossing}%
along a projection line (see \citet{Thiffeault2005,Thiffeault2010} and
Section~\ref{sec:braidfromdata}). %
\index{projection line|(}%
The \braid\ constructor allows us to do this easily.

\begin{figure}
\begin{center}
\subfigure[]{
  \includegraphics[height=.3\textheight]{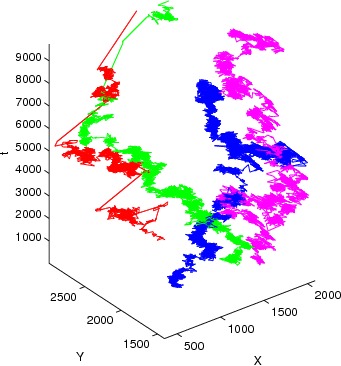}
  \label{fig:testdata_trajs3}
}\hspace{1em}
\subfigure[]{
  \includegraphics[height=.3\textheight]{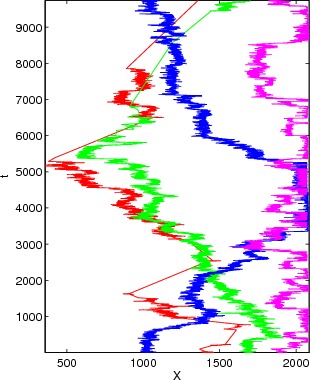}
  \label{fig:testdata_trajs}
}\hspace{1em}
\subfigure[]{
  \includegraphics[height=.3\textheight]{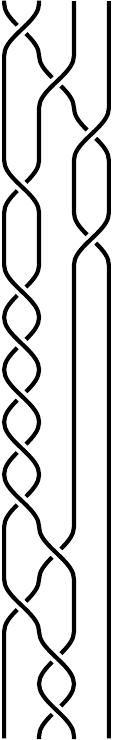}
  \label{fig:testdata_braid}
}\hspace{1em}
\subfigure[]{
  \includegraphics[height=.3\textheight]{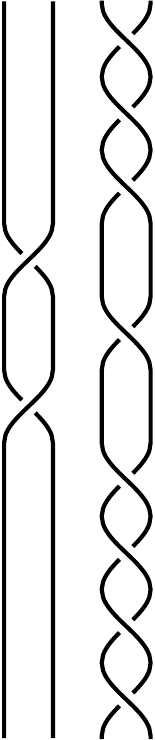}
  \label{fig:testdata_braidY}
}
\end{center}
\caption{(a) A dataset of four trajectories, (b) projected along the~$X$ axis.
  (c) The compacted braid~$ \sigma_2^{-2} \sigma_1^{-1} \sigma_2^{-1}
  \sigma_1^{-5} \sigma_3 \sigma_1^{-1}\sigma_3\sigma_2\sigma_1$ corresponding
  to the~$X$ projection in (b).  (d) The compacted
  braid~$\sigma_3^{-4}\sigma_1\sigma_3^{-1}\sigma_1\sigma_3^{-3}$
  corresponding to the~$Y$ projection, with closure enforced. %
  \index{braid!closure}%
  The braids in (c) and (d) are conjugate.\index{braid!conjugate}}
\end{figure}

The folder \lstinline{testsuite/testcases} %
\index{testsuite}%
contains a dataset of trajectories, from laboratory data for granular
media~\citep{Puckett2012}.  We load the
data:
\begin{lstlisting}[frame=single,framerule=0pt]
>> clear; load testdata
>> whos
  Name         Size               Bytes  Class     Attributes

  XY        9740x2x4             623360  double
  ti           1x9740             77920  double
\end{lstlisting}
Here \lstinline{ti} is the vector of times, and \lstinline{XY} is a
three-dimensional array: its first component specifies the timestep,
its second specifies the $X$ or $Y$ coordinate, and its third
specifies one of the~$4$ particles.  Figure~\ref{fig:testdata_trajs3}
shows
the~$X$ and~$Y$ coordinates of these four trajectories, with time
plotted vertically.  Figure~\ref{fig:testdata_trajs} shows the same
data, but projected along the~$X$ direction.  To construct a braid
from this data, we simply execute %
\index{braid class@\braid\ class!length@\lstinline{length}}
\begin{lstlisting}[frame=single,framerule=0pt]
>> warning('off','BRAIDLAB:braid:colorbraiding:notclosed')
>> b = braid(XY);
>> b.length

ans = 894
\end{lstlisting}
This is a very long braid!  (We've temporarily turned off a warning about the
data not being closed; we'll see what it's about in
Section~\ref{sec:projection}.)  But Figure~\ref{fig:testdata_trajs} suggests
that this is misleading: many of the crossings %
\index{crossing}%
are `wiggles' that cancel each other out.  Indeed, if we attempt to shorten
the braid: %
\index{braid class@\braid\ class!compact@\lstinline{compact}}%
\begin{lstlisting}[frame=single,framerule=0pt]
>> b = compact(b)

b = < -2 -2 -1 -2 -1 -1 -1 -1 -1  3 -1  3  2  1 >

>> b.length

ans = 14
\end{lstlisting}
we find the number of generators (the length) has dropped to~$14$!  We can
then plot this shortened braid as a braid diagram using \lstinline{plot(b)} %
\index{braid class@\braid\ class!plot@\lstinline{plot}} to produce
Figure~\ref{fig:testdata_braid}.  The braid diagram allows us to see some
topological information clearly, such as the fact that the second and third
particles undergo a large number of twists around each other; we can check
this by creating a subbraid %
\index{braid class@\braid\ class!subbraid class@\lstinline{subbraid}}%
with only those two strings:
\begin{lstlisting}[frame=single,framerule=0pt]
>> subbraid(b,[2 3])

ans = < -1 -1 -1 -1 -1 -1 -1 -1 >
\end{lstlisting}
which shows that the winding number between these two strings is~$-4$.

\subsubsection{Changing the projection line and enforcing closure}
\label{sec:projection}

The braid in the previous section was constructed from the data by assuming a
projection along the~$X$ axis (the default).  We can choose a different
projection by specifying an optional angle for the projection line; for
instance, to project along the~$Y$ axis we invoke
\begin{lstlisting}[frame=single,framerule=0pt]
>> b = braid(XY,pi/2);   % project onto Y axis
>> b.length

ans = 673

>> b.compact

ans = < -3 -3 -3 -3  1 -3 -3 -3 -3 >
\end{lstlisting}
In general, a change of projection line only changes the braid by
conjugation~\citep{Boyland1994,Thiffeault2010}.  We can test for
conjugacy: %
\index{braid!conjugate|(}%
\index{braid class@\braid\ class!compact@\lstinline{compact}|(}%
\index{braid class@\braid\ class!conjtest@\lstinline{conjtest}|(}%
\begin{lstlisting}[frame=single,framerule=0pt]
>> bX = compact(braid(XY,0)); bY = compact(braid(XY,pi/2));
>> conjtest(bX,bY)   % test for conjugacy of braids

ans = 0
\end{lstlisting}
\index{braid class@\braid\ class!compact@\lstinline{compact}|)}%
The braids are not conjugate.  This is because our trajectories do not form a
`true' braid: the final points do not correspond exactly with the initial
points, as a set.  This is the reason why we turned off a warning in
Section~\ref{sec:braidfromdataex}: \braidlab\ warns us when we're trying to
create a braid from data that doesn't `join up.'  The class
\lstinline{databraid} described in Section~\ref{sec:databraid} does not issue
this warning, since it is meant for noisy real-world data.

If we truly want a rotationally-conjugate braid out of our data, we need to
enforce a closure method: %
\index{braid!closure|(}%
\index{braid class@\braid\ class!closure@\lstinline{closure}|(}%
\index{braid class@\braid\ class!compact@\lstinline{compact}}%
\index{crossing!in braid closure|(}
\begin{lstlisting}[frame=single,framerule=0pt]
>> XY = closure(XY);   % close braid and avoid new crossings
>> bX = compact(braid(XY,0)), bY = compact(braid(XY,pi/2))

bX = < -2 -2 -1 -2 -1 -1 -1 -1 -1  3 -1  3  2  1 >

bY = < -3 -3 -3 -3  1 -3  1 -3 -3 -3 >
\end{lstlisting}
This default closure simply draws line segments from the final points to the
initial points in such a way that no new crossings are created in the~$X$
projection. %
\index{crossing!in braid closure|)}%
Hence, the $X$-projected braid \lstinline{bX} is unchanged by the closure, but
here the $Y$-projected braid \lstinline{bY} is longer by one generator
(\lstinline{bY} is plotted in Figure~\ref{fig:testdata_braidY}).  This is
enough to make the braids conjugate:
\begin{lstlisting}[frame=single,framerule=0pt]
>> [~,c] = conjtest(bX,bY)  % ~ means discard first return arg

c = < 3  2 >
\end{lstlisting}
\index{tilde@tilde (\lstinline{~}), as return argument}
where the optional second argument \lstinline{c} is the conjugating
braid, as we can verify:
\begin{lstlisting}[frame=single,framerule=0pt]
>> bX == c*bY*c^-1

ans = 1
\end{lstlisting}
There are other ways to enforce closure of a braid (see
\lstinline{help closure}), in particular
\lstinline{closure(XY,'MinDist')}, which minimizes the total distance
between the initial and final points.
\index{projection line|)}%
\index{braid!closure|)}%
\index{braid class@\braid\ class!closure@\lstinline{closure}|)}%

Note that \lstinline{conjtest} uses the library \emph{CBraid} \citep{CBraid} %
\index{CBraid}%
to first convert the braids to Garside canonical form \citep{Birman2005}, %
\index{braid!Garside form}%
then to determine conjugacy.  This is very inefficient, so is impractical for
large braids.
\index{braid!conjugate|)}%
\index{braid class@\braid\ class!conjtest@\lstinline{conjtest}|)}%

\subsubsection{The \lstinline{databraid} subclass}
\label{sec:databraid}

In some instances when dealing with data it is important to know the
\emph{crossing times}, %
\index{crossing!times|(}%
that is, the times at which two particles exchanged position along the
projection line. %
\index{projection line}%
A braid object does not keep this information, but there is an object that
does: a \lstinline{databraid}. %
\index{databraid class@\lstinline{databraid} class|(}%
Its constructor takes an optional vector of times as an argument, and it has a
data member \lstinline{tcross} %
\index{tcross@\lstinline{tcross}|(}%
that retains the crossing times.  Using the same data \lstinline{XY} from
before, sampled at times \lstinline{ti}, we have
\begin{lstlisting}[frame=single,framerule=0pt]
>> b = databraid(XY,ti);
>> b.tcross(1:3)

ans = 870.9010
      872.1758
      887.0089
\end{lstlisting}
\index{tcross@\lstinline{tcross}|)}%
Storing crossing times enables us to truncate generators in a
\lstinline{databraid} by retaining only those with crossing times within a
desired interval (see \lstinline{databraid.trunc}). There are always exactly
as many crossing times as generators in the braid. %
\index{databraid class@\lstinline{databraid} class!trunc@\lstinline{trunc}}%
\index{crossing!times|)}%

Many operations that can be done to a \lstinline{braid} also work on a
\lstinline{databraid}, with a few differences:
\begin{itemize}
\item\lstinline{compact} %
\index{braid class@\braid\ class!compact@\lstinline{compact}}%
works a bit differently.  It is less effective than \lstinline{braid.compact}
since it must preserve the order of generators in order to maintain the
ordering of the crossing times.
\item Equality testing checks if two \hbox{\lstinline{databraid}s} are
  lexicographically equal (i.e., generator-by-generator) and that their
  crossing times all agree.  This is very restrictive.  To check if the
  underlying braids are equal, first convert the \hbox{\lstinline{databraid}s}
  to \hbox{\lstinline{braid}s} by using the method
  \lstinline{databraid.braid}.
\item Multiplication of two \hbox{\lstinline{databraid}s} is only defined if
  the crossing times of the first braid are all earlier than the second.
\item Powers and inverses of \hbox{\lstinline{databraid}s} are not defined,
  since this would break the time-ordering of crossings.
\item%
  \index{Finite Time Braiding Exponent (FTBE)|(}%
  \index{braid class@\braid\ class!entropy@\lstinline{entropy}}%
  \index{braid class@\braid\ class!complexity@\lstinline{complexity}}%
  The entropy of a \lstinline{databraid} is an ambiguously defined concept.
  While entropy of certain braids can be computed non-iteratively, e.g., in
  \citet{Hall2009}, in general it is only estimated by an iterative process.
  Iterations rely on taking powers of the braid, which is not defined for
  \hbox{\lstinline{databraid}s}.  The functions \hbox{\lstinline{entropy}} and
  \hbox{\lstinline{complexity}} can still be used by converting
  \lstinline{databraid} objects to \lstinline{braid} objects; however, this
  should be avoided in favor of the appropriate concept for
  \hbox{\lstinline{databraid}s},
  the Finite Time Braiding Exponent (FTBE) \citep{Thiffeault2005,Budisic2015}.
  A \lstinline{databraid}~$b_T$ recorded over a time interval of length~$T$
  has an FTBE defined by
  \begin{equation}
    \FTBE(b_{T})
    = \frac{1}{T}\, \log \frac{ \lvert b_{T}\, \ell \rvert}{\lvert\ell\rvert},
    \label{eq:ftbe}
  \end{equation}
  where $\ell$ is a loop given by the generating set for the fundamental group
  of the disk with $n$ \index{punctures} punctures (see
  Section~\ref{sec:loopcoords}), and~$b_{T} \ell$ is that loop transformed by
  a single application of the braid.  Here~$\lvert\cdot\rvert$ is a measure of
  the length of the loop. Unless specified otherwise, \braidlab\ calculates
  \(T\) as the time elapsed between the first and last crossing in the
  braid. This duration could be much smaller than the length of trajectories
  analyzed, e.g., when no crossings occur near the beginning or the end of
  trajectories.  To set the custom value of \(T\) and other options, see
  documentation of the method \lstinline{databraid.ftbe}.
\end{itemize}

\index{databraid class@\lstinline{databraid} class|)}%
\index{databraid class@\lstinline{databraid} class!ftbe@\lstinline{ftbe}}%
\index{Finite Time Braiding Exponent (FTBE)|)}

\index{braid class@\braid\ class!constructor!from data|)}
\index{braid!from data|)}%
\index{braid class@\braid\ class|)}

\subsection{The \loopc\ class}
\label{sec:loop}

\index{loop class@\loopc\ class|(}

\subsubsection{Loop coordinates}
\label{sec:loopc}

\begin{figure}
\begin{center}
\subfigure[]{
  \includegraphics[height=.22\textheight]{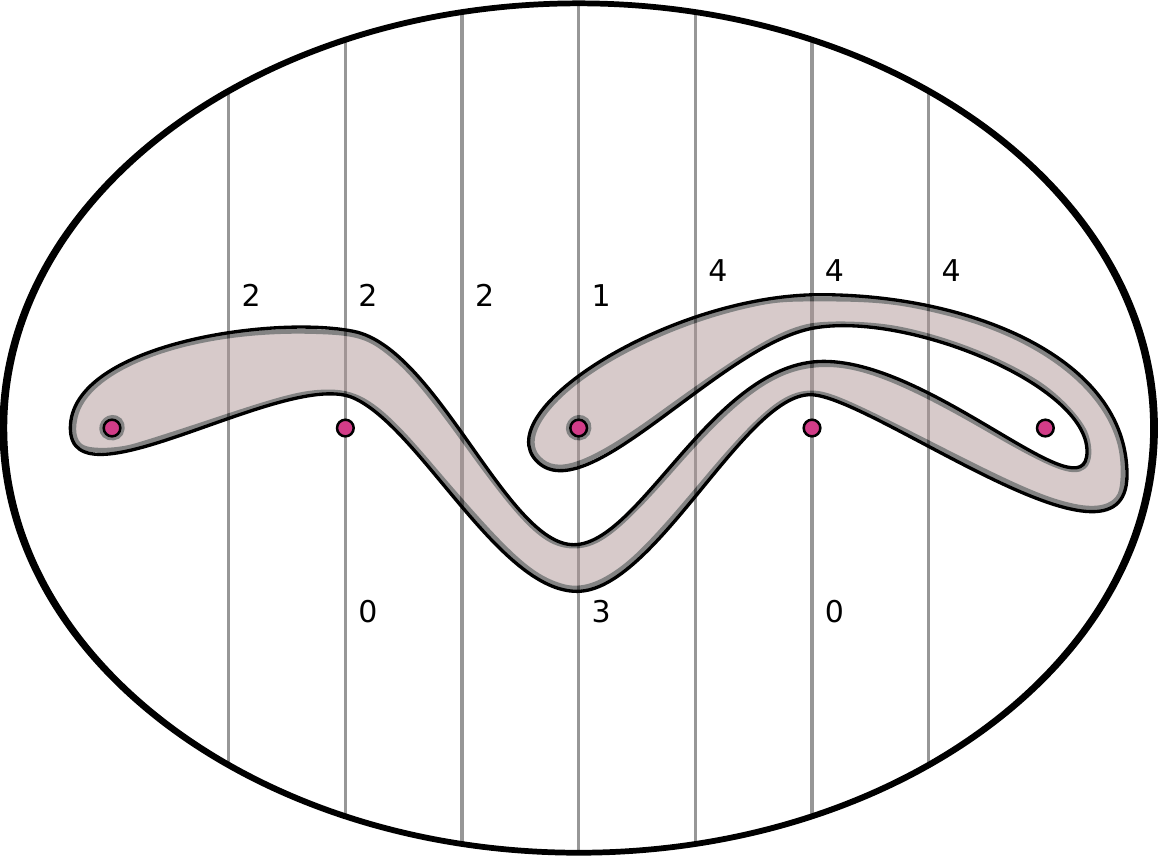}
  \label{fig:dynn_loop}
}\hspace{1em}
\subfigure[]{
  \includegraphics[height=.22\textheight]{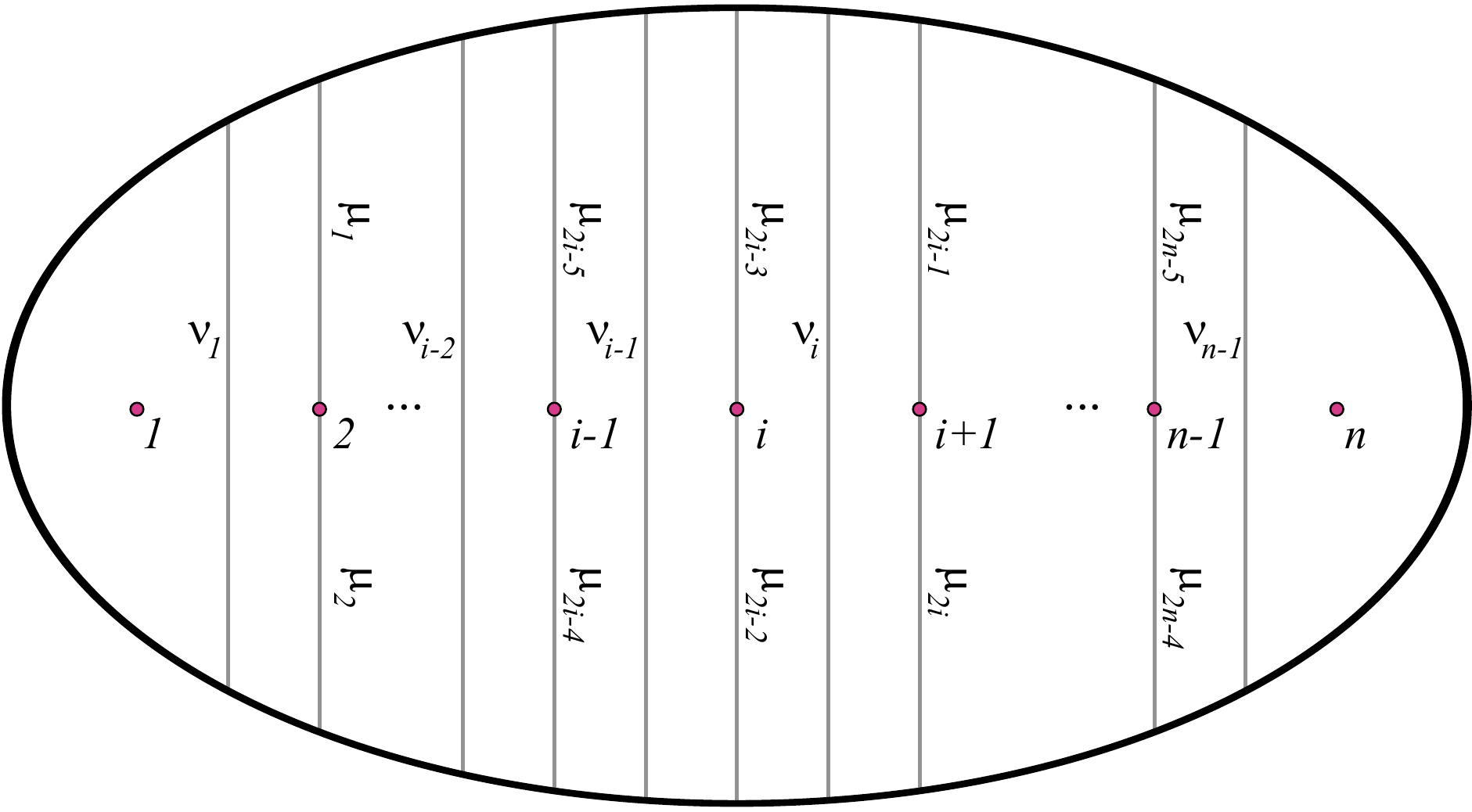}
  \label{fig:dynn_def}
}
\end{center}
\caption{(a) A simple close loop in a disk with~$\nn=5$ punctures.
  (b) Definition of intersection numbers~$\mu_i$ and~$\nu_i$.
  [From~\citet{Thiffeault2010}.] \index{loop!intersection numbers}}
\end{figure}

A simple closed loop on a disk with~$5$ punctures %
\index{disk, punctured}%
\index{punctures}%
is shown in Figure~\ref{fig:dynn_loop}.
We consider equivalence classes of such loops under homotopies %
\index{loop!homotopy classes}%
relative to the punctures.  %
\index{punctures}%
In particular, the loops are \emph{essential}, %
\index{loop!essential}%
meaning that they are not null-homotopic or homotopic to the boundary or a
puncture.  The \emph{intersection numbers} %
\index{loop!intersection numbers|(}%
\index{intersection numbers|see{loop intersection numbers}}%
are also shown in Figure~\ref{fig:dynn_loop}: these count the minimum number
of intersections of an equivalence class of loops with the fixed vertical
lines shown.  For~$\nn$ punctures, we define the intersection numbers~$\mu_i$
and~$\nu_i$ in Figure~\ref{fig:dynn_def}.

Any given loop will lead to a unique set of intersection numbers, but
a general collection of intersection numbers do not typically
correspond to a loop. %
\index{loop!coordinates|(}
It is therefore more convenient to define
\begin{equation}
  \ac_\ip = \tfrac12\l(\mu_{2\ip} - \mu_{2\ip-1}\r), \qquad
  \bc_\ip = \tfrac12\l(\nu_\ip - \nu_{\ip+1}\r), \qquad
  \ip=1,\ldots,\nn-2.
\end{equation}
We then combine these in a vector of length~$(2\nn-4)$,
\begin{equation}
  \abv = (\ac_1,\ldots,\ac_{\nn-2},\bc_1,\ldots,\bc_{\nn-2}),
  \label{eq:abvdef}
\end{equation}
which gives the \emph{loop coordinates} (or \emph{Dynnikov coordinates}) for
the loop.  (Some authors such as~\citet{Dehornoy2008} give the coordinates
as~$(\ac_1,\bc_1,\ldots,\ac_{\nn-2},\bc_{\nn-2})$.)  There is now a bijection
between~$\mathbb{Z}^{2\nn-4}$ and essential simple closed
loops~\citep{Dynnikov2002,Moussafir2006,Hall2009,Thiffeault2010}.  Actually, %
\index{multiloop|see{loop, multi-}}%
\index{loop!multi-} \emph{multiloops}: loop coordinates can describe unions of
disjoint loops (see Section~\ref{sec:loopcoords}).%
\footnote{%
  Here we use multiloop \index{loop!multi-} as a convenient mnemonic.  The
  technical term is \emph{integral lamination}: %
  \index{integral lamination|see{loop, multi-}}%
  a set of disjoint non-homotopic simple closed curves~\citep{Moussafir2006}.}

\index{loop class@\loopc\ class!constructor|(}
Let's create the loop in Figure~\ref{fig:dynn_loop} as a \loopc\ object:
\begin{lstlisting}[frame=single,framerule=0pt]
>> l = loop([-1 1 -2 0 -1 0])

l = (( -1 1 -2 0 -1 0 ))
\end{lstlisting}
\index{loop class@\loopc\ class!constructor|)}
Figure~\ref{fig:dynn_loop2} shows the output of the \lstinline{plot(l)} %
\index{loop class@\loopc\ class!plot@\lstinline{plot}}
command.  We can convert from loop coordinates\index{loop!coordinates} to
intersection numbers with
\begin{lstlisting}[frame=single,framerule=0pt]
>> intersec(l)

ans = 2 0 1 3 4 0 2 2 4 4   % [mu1 ... mu6 nu1 ... nu4]
\end{lstlisting}
which returns~$\mu_1\dots\mu_{2n-4}$ followed by~$\nu_1\dots\mu_{n-1}$, as
defined in Figure~\ref{fig:dynn_def}.
\index{loop!intersection numbers|)}%

We can also extract the loop coordinates from a \loopc\ object using the
methods \lstinline{a}, \lstinline{b}, and \lstinline{ab}: %
\index{loop class@\loopc\ class!%
  abab@\lstinline{a}, \lstinline{b}, \lstinline{ab}}%
\begin{lstlisting}[frame=single,framerule=0pt]
>> l = loop ([-1 1 -2 0 -1 0]);
>> l.a

ans = -1     1    -2

>> l.b

ans =  0    -1     0

>> [a,b] = l.ab

a = -1     1    -2
b =  0    -1     0
\end{lstlisting}
As for braids, \lstinline{l.n} returns the number of punctures (or strings). %
\index{loop class@\loopc\ class!number of punctures (\lstinline{n})}

\subsubsection{Acting on loops with braids}
\label{sec:actingonloops}

Now we can act on this loop with braids. %
\index{braid class@\braid\ class!action on \loopc\ (\lstinline{*})|(}%
For example, we define the braid
\lstinline{b} to be~$\sigma_1^{-1}$ with~$5$ strings, corresponding to the~$5$
punctures, %
\index{punctures}%
\index{braid class@\braid\ class!multiplication (\lstinline{*})|(}
and then act on the loop \lstinline{l} by using the multiplication operator:
\begin{figure}
\begin{center}
\subfigure[]{
  \includegraphics[width=.6\textwidth]{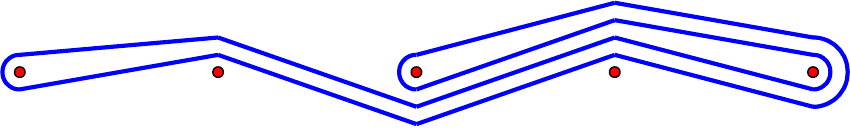}
  \label{fig:dynn_loop2}
}\hspace{1em}
\subfigure[]{
  \includegraphics[width=.6\textwidth]{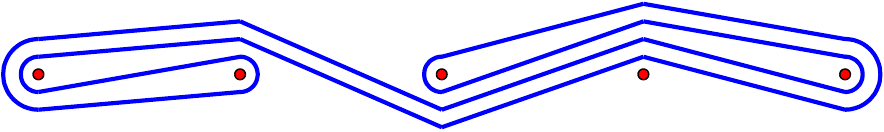}
  \label{fig:dynn_loop2_sigm1}
}
\end{center}
\caption{(a) The loop \lstinline{((-1 1 -2 0 -1 0))}.  (b) The braid generator
  $\sigma_1^{-1}$ applied to the loop in (a).}
\end{figure}
\begin{lstlisting}[frame=single,framerule=0pt]
>> b = braid([-1],5);   % one generator with 5 strings
>> b*l                  % act on a loop with a braid

ans = (( -1  1 -2  1 -1  0 ))
\end{lstlisting}
Figure~\ref{fig:dynn_loop2_sigm1} shows \lstinline{plot(b*l)}.  The first and
second punctures %
\index{punctures|(}%
were interchanged counterclockwise (the action of~$\sigma_1^{-1}$), dragging
the loop along. %
\index{loop!coordinates|)}%
\index{braid class@\braid\ class!multiplication (\lstinline{*})|)}%
\index{braid class@\braid\ class!action on \loopc\ (\lstinline{*})|)}%

\index{loop class@\loopc\ class!minlength@\lstinline{minlength}|(}
\index{loop!minimum length} The minimum length of an equivalence class of
loops is determined by assuming the punctures are one unit of length apart and
have zero size.  After pulling tight the loop on the punctures, %
\index{punctures|)}%
it is then made up of unit-length segments.  The minimum length is thus an
integer.  For the loop in Figure~\ref{fig:dynn_loop2},
\begin{lstlisting}[frame=single,framerule=0pt]
>> minlength(l)

ans = 12
\end{lstlisting}
\index{loop class@\loopc\ class!minlength@\lstinline{minlength}|)}
\index{loop class@\loopc\ class!intaxis@\lstinline{intaxis}|(}
Another useful measure of a loop's complexity is its minimum intersection
number with the real axis~\citep{Moussafir2006,Hall2009,Thiffeault2010}, which
for this loop is the same as its minimum length:
\begin{lstlisting}[frame=single,framerule=0pt]
>> intaxis(l)

ans = 12
\end{lstlisting}
The \lstinline{intaxis} method is used to measure a braid's geometric
complexity, %
\index{braid!complexity|(}%
\index{braid class@\braid\ class!complexity@\lstinline{complexity}|(}%
as defined by~\citet{Dynnikov2007}.
\index{loop class@\loopc\ class!intaxis@\lstinline{intaxis}|)}

\index{loop class@\loopc\ class!constructor|(}
\index{loop class@\loopc\ class!vectorized|(}
Sometimes we wish to study a large set of different loops.  The loop
constructor vectorizes:
\begin{lstlisting}[frame=single,framerule=0pt]
>> ll = loop([-1 1 -2 0; 1 -2 3 4])

ll = (( -1  1 -2  0 ))
     ((  1 -2  3  4 ))
\end{lstlisting}
\index{loop class@\loopc\ class!constructor|)}
\index{loop class@\loopc\ class!minlength@\lstinline{minlength}|(}
We can then, for instance, compute the length of every loop:
\begin{lstlisting}[frame=single,framerule=0pt]
>> minlength(ll)

ans = 14
      34
\end{lstlisting}
\index{loop class@\loopc\ class!minlength@\lstinline{minlength}|)}
or even act on all the loops with the same braid:
\index{braid class@\braid\ class!action on \loopc\ (\lstinline{*})|(}%
\index{braid class@\braid\ class!multiplication (\lstinline{*})|(}
\begin{lstlisting}[frame=single,framerule=0pt]
>> b = braid([1 -2]);
>> b*ll

ans = (( 2   1  -2   1 ))
      (( 5  -2  -3  11 ))
\end{lstlisting}
\index{braid class@\braid\ class!multiplication (\lstinline{*})|)}
\index{braid class@\braid\ class!action on \loopc\ (\lstinline{*})|)}%
Some commands, such as \lstinline{plot}, do not vectorize.  Different loops
can then be accessed by indexing, such as~\lstinline{plot(ll(2))}.
\index{loop class@\loopc\ class!vectorized|)}

The \lstinline{entropy} method %
\index{braid class@\braid\ class!entropy@\lstinline{entropy}|(}%
\index{braid!entropy|(}%
of the \lstinline{braid} class
(Section~\ref{sec:braidclass}) computes the topological entropy of a braid by
repeatedly acting on a loop, and monitoring the growth rate of the loop.  For
example, let us compare the entropy obtained by acting~$100$ times on an
initial loop, compared with the \lstinline{entropy} method:
\begin{lstlisting}[frame=single,framerule=0pt]
>> b = braid([1 2 3 -4]);
% apply braid 100 times to l, then compute growth of length
>> log(minlength(b^100*l)/minlength(l)) / 100

ans = 0.7637

>> entropy(b)

ans = 0.7672
\end{lstlisting}
The entropy value returned by \lstinline{entropy(b)} is more precise,
since that method monitors convergence and adjusts the number of
iterations accordingly. %
\index{braid class@\braid\ class!entropy@\lstinline{entropy}|)}%
\index{braid!entropy|)}%

\subsection{Loop coordinates for a braid}
\label{sec:loopcoords}

\index{braid!loop coordinates|(}
\index{loop!coordinates|(}
\index{loop class@\loopc\ class!constructor|(}

\begin{figure}
\begin{center}
\subfigure[]{
  \includegraphics[width=.7\textwidth]{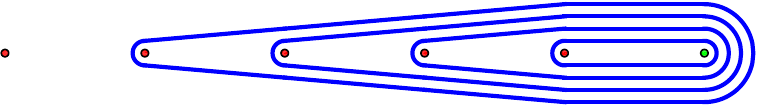}
  \label{fig:fundloops}
}\hspace{1em}
\subfigure[]{
  \includegraphics[width=\textwidth]{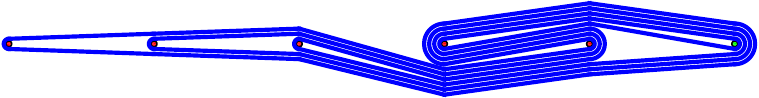}
  \label{fig:fundloops_act}
}
\end{center}
\caption{(a) The multiloop created by \lstinline{loop(5,'BP')}, with basepoint
  puncture in green.  (b) The multiloop \lstinline{b*loop(5,'BP')}, where
  \lstinline{b} is the braid
  $\sigma_1\sigma_2\sigma_3\sigma_4^{-1}$. \index{loop!multi-}}
\end{figure}

The command \lstinline{loop(n,'BasePoint')} returns a \emph{canonical set of
  loops} for~$n$ punctures:%
\index{punctures|(}%
\begin{lstlisting}[frame=single,framerule=0pt]
>> l = loop(5,'BP')    % 'BP' is short for 'BasePoint'

ans = (( 0  0  0  0 -1 -1 -1 -1 ))*
\end{lstlisting}
\index{loop!multi-|(} This multiloop is depicted in
Figure~\ref{fig:fundloops}, with basepoint puncture shown in green.  The
\lstinline{*} indicates that this loop has a basepoint.  Note that the
multiloop returned by \hbox{\lstinline{loop(5,'BP')}} actually has 6
punctures!  The rightmost puncture is meant to represent the boundary of a
disk, %
\index{disk, punctured}%
or a base point for the fundamental group on a sphere with $n$ punctures.  The
loops form a (nonoriented) generating set for the fundamental group of the
disk with $n$ punctures.  The extra puncture thus plays no role dynamically,
and \lstinline{l.n} returns~5.  If you want the true total number of
punctures, including the base point, use \lstinline{l.totaln}. %
\index{loop class@\loopc\ class!totaln@\lstinline{totaln}}
\index{loop class@\loopc\ class!constructor|)}
\index{punctures|)}%
\index{loop!multi-|)}

The canonical set of loops allows us to define loop coordinates for a braid,
which is a unique normal form.
The canonical loop coordinates for braids exploit the fact that two braids are
equal if and only if they act the same way on the fundamental group of the
disk \citep{Dehornoy2008}.  Hence, if we take a braid and act on
\lstinline{loop(5,'BP')},
\begin{lstlisting}[frame=single,framerule=0pt]
>> b = braid([1 2 3 -4]);
>> b*loop(5,'BP')

ans = (( 0  0  3 -1 -1 -1 -4  3 ))*
\end{lstlisting}
then the set of numbers \lstinline{(( 0 0 3 -1 -1 -1 -4 3 ))*} can be thought
of as \emph{uniquely} characterizing the braid.  It is this property that is
used to rapidly determine equality of braids.  (The loop
\lstinline{b*loop(5,'BP')} is plotted in Figure~\ref{fig:fundloops_act}.)  The
same loop coordinates for the braid can be obtained without creating an
intermediate loop with %
\index{braid class@\braid\ class!loopcoords@\lstinline{loopcoords}}
\begin{lstlisting}[frame=single,framerule=0pt]
>> loopcoords(b)

ans = (( 0  0  3 -1 -1 -1 -4  3 ))*
\end{lstlisting}

\index{loop class@\loopc\ class|)}
\index{loop!coordinates|)}
\index{braid!loop coordinates|)}

\jlt{Next section: Braid from random walks?  Compute runs of same gen.}

\section{The effective linear action and its cycles}
\label{sec:elacycles}

\subsection{Effective linear action}
\label{sec:linact}

\index{effective linear action|(}%
\index{linear action|see{effective linear action}}%
\index{action!effective linear|(}%
\index{action!of braid on loop|(}%

In Section~\ref{sec:actingonloops} we introduced the action of a
braid~$\gamma$ on a loop~$\abv$.
Here~$\abv=(\ac_1,\ldots,\ac_{\nn-2},\bc_1,\ldots,\bc_{\nn-2})$ is a vector of
coordinates for the loop, %
\index{loop!coordinates}%
defined in Section~\ref{sec:loopc}.  We write~$\abv' = \gamma\cdot\abv$ for
the new, updated coordinates after the action.  These updated coordinates are
given by composing the action of individual generators.

\index{update rules|see{action}}
\index{action!update rules|(}
For~$1 < \ip < \nn-1$, we can express the update rules for the braid group
generator~$\sigma_\ip$ acting on~$\abv$ as
\begin{subequations}
\begin{align}
  \acnew_{\ip-1} &= \ac_{\ip-1} - \pos{\bc_{\ip-1}}
    - \pos{\l(\pos{\bc_\ip} + \cc_{\ip-1}\r)}\,,\\
  \bcnew_{\ip-1} &= \bc_\ip + \neg{\cc_{\ip-1}}\,,\\
  \acnew_\ip &= \ac_\ip - \neg{\bc_\ip}
    - \neg{\l(\neg{\bc_{\ip-1}} - \cc_{\ip-1}\r)}\,,\\
  \bcnew_\ip &= \bc_{\ip-1} - \neg{\cc_{\ip-1}}\,,
\end{align}
\label{eq:ur}%
\end{subequations}
where
\begin{equation}
  \cc_{\ip-1} = \ac_{\ip-1} - \ac_\ip - \pos{\bc_\ip} + \neg{\bc_{\ip-1}}\,.
  \label{eq:ccdef}
\end{equation}
Coordinates not listed (i.e., $\ac_k$ and~$\bc_k$ for~$k\ne\ip$ or~$\ip-1$)
are unchanged.  The superscripts~${}^{+/-}$ are defined as
\begin{equation}
  \pos\fc \ldef \max(\fc,0),\qquad
  \neg\fc \ldef \min(\fc,0).
\end{equation}
(See~\citet{Thiffeault2010} for the update rules for the
generators~$\sigma_1$, $\sigma_{\nn-1}$, and the inverse generators.  The
update rules are in several other papers but use different conventions.)
\index{action!update rules|)}

Notice that the action~\eqref{eq:ur} is \emph{piecewise-linear} in the loop
coordinates: once the~${}^{+/-}$ operators are resolved, what is left is a
linear operation on the vector~$\abv$.  We can thus write
\begin{equation}
  \abv' = M(\gamma,\abv)\cdot\abv,
  \qquad
  M(\gamma,\abv) \in \mathrm{SL}_{2\nn-4}(\mathbb{Z}),
\end{equation}
where the dot now denotes the standard matrix product.  Here~$M(\gamma,\abv)$
is the \emph{effective linear action} of the braid~$\gamma$ on the
loop~$\abv$.

\index{braid class@\braid\ class!action on \loopc\ (\lstinline{*})|(}%
Let's show an example using \braidlab.  We take the
braid~$\sigma_1\sigma_2^{-1}$ and the loop with coordinates~$\ac_1=0$,
$\bc_1=-1$.  The action is
\begin{lstlisting}[frame=single,framerule=0pt]
>> b = braid([1 -2]); l = loop([0 -1]);
>> lp = b*l

lp = (( 1 -1 ))
\end{lstlisting}
The effective linear action can be obtained by requesting a second output
argument from the result of~\lstinline{*}:
\begin{lstlisting}[frame=single,framerule=0pt]
>> [lp,M] = b*l; full(M)

ans = 1 -1
      0  1
\end{lstlisting}
Note that the effective linear action \lstinline{M} is by default returned as
a sparse matrix, %
\index{sparse matrix}%
which it often is when dealing with many strands.  We use \lstinline{full} %
\index{full@\lstinline{full}}%
to convert it back into a regular full matrix.  We can then verify that the
matrix product of \lstinline{M} and the column vector of coordinates
\lstinline{l.coords'} is the same as the action \lstinline{lp = b*l}:
\begin{lstlisting}[frame=single,framerule=0pt]
>> M*l.coords'

ans = 1
     -1

>> lp.coords'

ans = 1
     -1
\end{lstlisting}
The difference is that \lstinline{M} may only be applied \emph{to this
  specific loop} (or a loop that happens to share the same effective linear
action).

A common thing to do is to find the effective linear action on the canonical
set \lstinline{loop(b.n,'BP')} (see Section~\ref{sec:loopcoords}):
\begin{lstlisting}[frame=single,framerule=0pt]
>> [~,M] = b*loop(b.n,'BP'); full(M)

ans = 0     0    -1     0
      0     1     0     1
      0     1     1     1
      1    -1    -1     0
\end{lstlisting}
The canonical set assumes an extra puncture, so the matrix dimension is larger
by~$2$.

The effective linear action doesn't seem to offer much at this point.  Its
real advantage will become apparent in Section~\ref{sec:cycles}, when we find
that it can achieve periodic limit cycles.

\index{braid class@\braid\ class!action on \loopc\ (\lstinline{*})|)}%
\index{action!of braid on loop|)}%

\subsection{Limit cycles of the effective linear action}
\label{sec:cycles}

\index{cycle!of effective linear action|(}
\index{effective linear action!cycle|(}

\begin{figure}
\begin{center}
\subfigure[]{
  \includegraphics[height=.35\textheight]{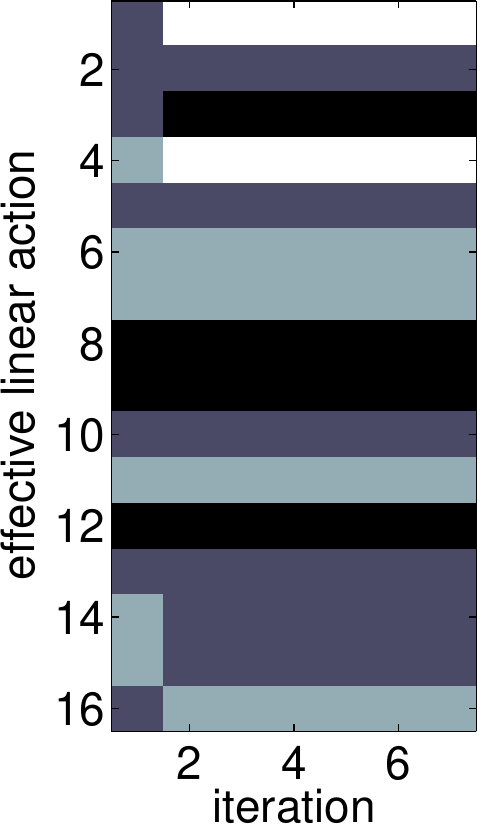}
  \label{fig:efflinact1}
}\hspace{.5em}
\subfigure[]{
  \includegraphics[height=.35\textheight]{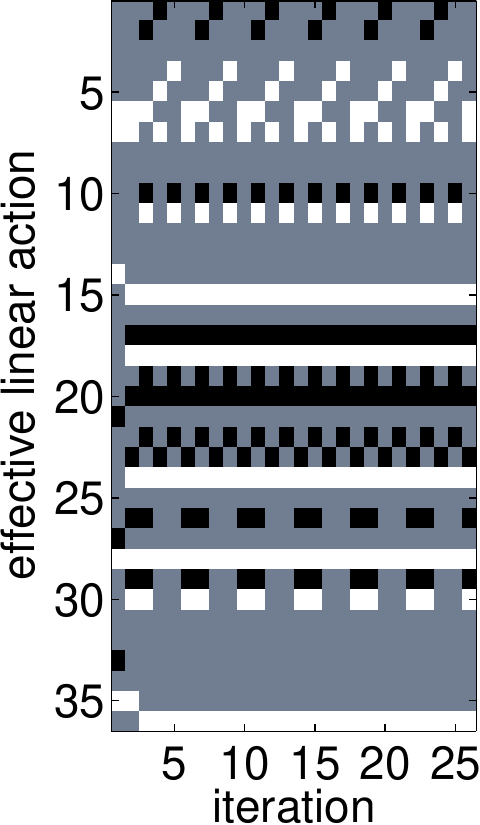}
  \label{fig:efflinact2}
}\hspace{.5em}
\subfigure[]{
  \includegraphics[height=.35\textheight]{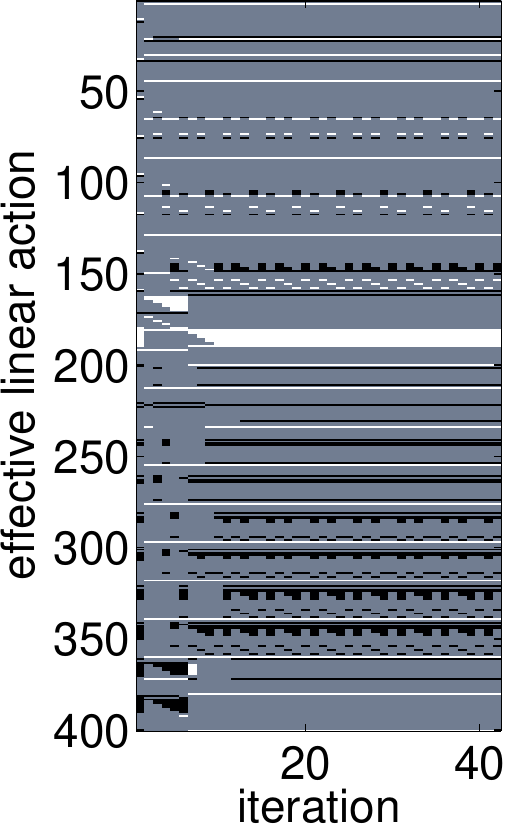}
  \label{fig:efflinact3}
}
\end{center}
\caption{The plot produced by \lstinline{cycle(b,'Plot')} %
  for (a) \lstinline{b = braid([1 -2])}; (b) \lstinline{b = braid([1 2 3])};
  (c) \lstinline{b = braid('Psi',11)}.}
\index{braid class@\braid\ class!cycle@\lstinline{cycle}}%
\label{fig:efflinact}
\end{figure}

The effective linear action has a very interesting behavior when a braid is
iterated on some initial loop.  Consider the following example:
\begin{lstlisting}[frame=single,framerule=0pt]
>> b = braid([1 -2]); l = loop([1 1]);
>> [l,M] = b*l; l, full(M)

l = (( 3 -1 ))

M = 2   1
   -1   0
\end{lstlisting}
Now repeat this last command:
\begin{lstlisting}[frame=single,framerule=0pt]
>> [l,M] = b*l; l, full(M)

l = (( 7 -4 ))

M = 2  -1
   -1   1
\end{lstlisting}
And again:
\begin{lstlisting}[frame=single,framerule=0pt]
>> b = braid([1 -2]); l = loop([1 1]);
>> [l,M] = b*l; l, full(M)

l = (( 18 -11 ))

M = 2  -1
   -1   1
\end{lstlisting}
The effective linear action \lstinline{M} has not changed.  In fact it has
achieved a fixed point: %
\index{fixed point|see{cycle}}%
running the same command again will change the loop, but the linear action
will remain the same forever.  \braidlab\ can automate the iteration with the
method \lstinline{cycle}. %
\index{braid class@\braid\ class!cycle@\lstinline{cycle}|(}%
Figure~\ref{fig:efflinact1} shows the output of
\begin{lstlisting}[frame=single,framerule=0pt]
>> b = braid([1 -2]); M = cycle(b,'Plot');
\end{lstlisting}
The member function \lstinline{cycle} iterates the braid on an initial loop,
taken to be the canonical set \lstinline{loop(b.n,'BP')}.  The vertical axis
in Fig.~\ref{fig:efflinact1} shows the elements of the effective linear action
as a function of iterates of the braid.  The matrix of the action is flattened
into a vector of length~$4^2$, where~$4$ is the dimension the initial loop
\lstinline{loop(b.n,'BP')}.  It is evident that the fixed point is reached
rapidly, since the `stripes' stop changing.

\index{braid!pseudo-Anosov|(}%
Such fixed points of the effective linear action are ubiquitous for braids
corresponding to a pseudo-Anosov isotopy class, such
as~$\sigma_1\sigma_2^{-1}$. %
In general, instead of a fixed point we may find a \emph{limit cycle} of some
period.  \citet{Yurttas2014_preprint} discussed these limit cycles for
pseudo-Anosov braids: they occur when the unstable foliation falls on the
boundary of the linear regions of the update rules.  We can reproduce her
example with the following:%
\footnote{To get exactly the same matrices, we use the
  braid~$\sigma_1^{-1}\sigma_2^{-1}\sigma_3^{-1}\sigma_4$ rather than
  her~$\sigma_1\sigma_2\sigma_3\sigma_4^{-1}$, since her generators rotate the
  punctures counterclockwise.}
\begin{lstlisting}[frame=single,framerule=0pt]
>> b = braid([-1 -2 -3 4]);
>> M = b.cycle(loop(b.n),'Iter')

M = [6x6 double]    [6x6 double]
\end{lstlisting}
The option \lstinline{'Iter'} tells \lstinline{cycle} to compute an individual
matrix for each iterate of the cycle, rather than the net product of all the
matrices in the cycle.  The output is a cell array %
\index{Matlab!cell array}%
of two~$6$ by~$6$ matrices, corresponding to the period-2 cycle:
\begin{lstlisting}[frame=single,framerule=0pt]
>> full(M{1}), full(M{2})

ans = -1     1     0     0     0     0
       0     0     0     1     1     0
       0     0     2    -1    -1     1
       0     0     0     0     1     0
      -1     0     1    -1    -1     1
       0     0     1     0     0     1

ans =  0     0     0     1     0     0
       0     0     0     1     1     0
       0     0     2    -1    -1     1
      -1     1     0    -1     1     0
       0    -1     1     0    -1     1
       0     0     1     0     0     1
\end{lstlisting}
as given by \citet{Yurttas2014_preprint}.  Note that we use an initial loop
for~$\nn$ punctures (\lstinline{loop(b.n)}) without base point, rather than
the default, to reproduce her example exactly.  For the pseudo-Anosov case,
any initial loop will give the same matrices.

What is more surprising is that these limit cycles occur for finite-order
braids as well.  Figure~\ref{fig:efflinact2} is produced by
\begin{lstlisting}[frame=single,framerule=0pt]
>> b = braid([1 2 3]); [~,period] = cycle(b,'Plot')

period = 4
\end{lstlisting}
Indeed, staring at the pattern in Fig.~\ref{fig:efflinact2} it is easy to see
that the effective action does achieve a limit cycle of period~$4$.  This
braid is definitely not pseudo-Anosov: it is finite-order.
\index{braid!finite-order}%
However, we do not expect such limit cycles to be unique in the
non-pseudo-Anosov case.

Pseudo-Anosov braids %
\index{braid!pseudo-Anosov}%
can achieve longer cycles, which \braidlab\ can find:
Figure~\ref{fig:efflinact3} is the plot produced by
\begin{lstlisting}[frame=single,framerule=0pt]
>> b = braid('Psi',11); [M,period] = cycle(b,'Plot');
\end{lstlisting}
The period here is~$5$, and the matrix~$M$ is~$20$ by~$20$.  The braid %
\index{psi braids@$\psi$ braids}%
\index{braid!entropy!minimum}%
\lstinline{braid('Psi',11)} is the braid~$\psi_{11}$ in the notation of
\citet{Venzke_thesis}.  It is a pseudo-Anosov braid with low %
\index{braid!pseudo-Anosov}%
\index{dilatation}%
dilatation~\citep{Hironaka2006,Thiffeault2006}, conjectured to be the lowest
possible for~$11$ strings.%
\footnote{The dilatation of a braid is the exponential of its entropy.}
The braids~$\psi_\nn$ are known to have to lowest dilatation for~$\nn$ string
for $\nn\le8$~\citep{LanneauThiffeault2011_braids}.

The largest eigenvalue of the matrix \lstinline{M} gives us the
\index{dilatation}%
dilatation of the braid, which in itself is not a real improvement over our
earlier entropy iterative algorithm (Section~\ref{sec:entropy}).  However,
with the matrix in hand we can find the characteristic polynomial:%
\footnote{Matlab's symbolic toolbox %
\index{Matlab!symbolic toolbox}%
is needed for \lstinline{poly2sym} and \lstinline{factor}.}
\begin{lstlisting}[frame=single,framerule=0pt]
>> b = braid('Psi',7); [M,period] = cycle(b);
>> factor(poly2sym(charpoly(M)))  % convert to symbolic form

ans = (x^2 + 1)*(x^3 - x^2 - 1)*(x^3 + x - 1)*(x - 1)^2*(x + 1)^2
\end{lstlisting}
\index{braid class@\braid\ class!cycle@\lstinline{cycle}|)}%
Compare this to the known polynomial that gives the dilation:
\index{psiroots@\lstinline{psiroots}}
\begin{lstlisting}[frame=single,framerule=0pt]
>> factor(poly2sym(psiroots(7,'Poly')))

ans = (x + 1)*(x^3 - x^2 - 1)*(x^3 + x - 1)
\end{lstlisting}
(The function \lstinline{psiroots} returns the roots and characteristic
polynomial of a~$\psi$ braid; this is useful for testing purposes.)  Note that
the factor whose largest root is the dilatation, %
\index{dilatation}%
\lstinline{x^3 - x^2 - 1}, appears in both polynomials.  This is not always
the case, though the dilatation has to be a root of both polynomials.
\index{braid!pseudo-Anosov|)}%

To our knowledge, the existence of these limit cycles has not been fully
explained (except in the pseudo-Anosov case by \citet{Yurttas2014_preprint}).
They seem to occur for \emph{any} braid, regardless of its %
\index{homeomorphism!isotopy classes}%
isotopy class.  In that sense they could provide an alternative to the the
Bestvina--Handel train track algorithm~\citep{Bestvina1995}, %
\index{Bestvina--Handel algorithm}%
which is used to compute the isotopy class of a braid.

\index{cycle!of effective linear action|)}%
\index{effective linear action|)}%
\index{effective linear action!cycle|)}
\index{action!effective linear|)}%
\index{limit cycle|see{cycle}}

\section{An example: Taffy pullers}
\label{sec:taffy}

\index{taffy pullers|(}

\begin{figure}
\begin{center}
\subfigure[]{
  \includegraphics[height=.2\textheight]{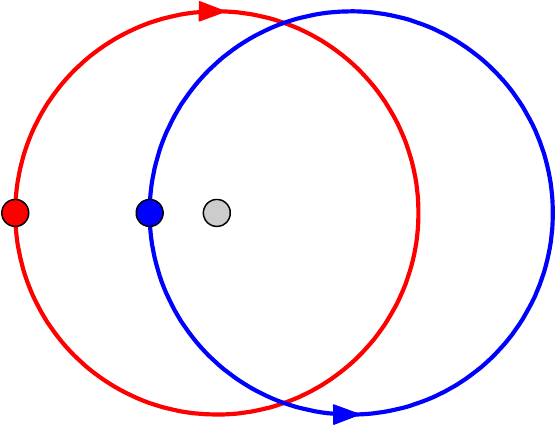}
  \label{fig:taffy_3rods}
}\hspace{1em}
\subfigure[]{
  \includegraphics[height=.2\textheight]{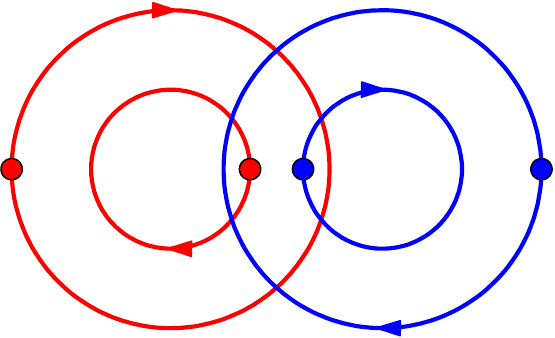}
  \label{fig:taffy_4rods}
}
\end{center}
\caption{(a) Three-rod taffy puller.  (b) Four-rod taffy puller.}
\label{fig:taffy}
\end{figure}

Taffy pullers are a class of devices designed to stretch and fold soft candy
repeatedly \citep{MattFinn2011_silver,Thiffeault2018}.  The goal is to aerate
the taffy.  Since many folds are required, the process has been mechanized
using fixed and moving rods.  The two most typical designs are shown in
Figure~\ref{fig:taffy}: the one in
Figure~\ref{fig:taffy_3rods} has a single fixed rod (gray) and two moving
rods, each rotating on a different axis.  The design in
Figure~\ref{fig:taffy_3rods} has four moving rods, sharing two axes of
rotation.  (There are several videos of taffy pullers on
\href{http://www.youtube.com/watch?v=6QkGp2qBbn4}{YouTube}.)

Let's use \braidlab\ to analyze the rod motion.  From the folder
\lstinline{doc/examples}, run the command%
\footnote{When using the parallel code, the generator sequences in this
  section may sometimes differ from run-to-run, due to the simultaneous
  crossings.  However, the braids themselves are still equal, after the
  relations~\eqref{eq:relations} are taken into account.}
\index{taffy@\lstinline{taffy}|(}%
\begin{lstlisting}[frame=single,framerule=0pt]
>> b = taffy('3rods')

b = < -2  1  1 -2 >
\end{lstlisting}
which also produces Figure~\ref{fig:taffy_3rods}.  The Thurston--Nielsen %
\index{braid!Thurston--Nielsen type \lstinline{tntype}}%
type and topological entropy %
\index{braid!entropy}%
of this braid are %
\index{braid class@\braid\ class!train@\lstinline{train}|(}%
\begin{lstlisting}[frame=single,framerule=0pt]
>> train(b)

ans = struct with fields:

       braid: [1x1 braidlab.braid]
      tntype: 'pseudo-Anosov'
     entropy: 1.7627
    transmat: [2x2 double]
       ttmap: {5x1 cell}
\end{lstlisting}
One would expect a competent taffy puller to be pseudo-Anosov, %
\index{braid!pseudo-Anosov}%
as this one is.  It implies that there is no `bad' initial condition where a
piece of taffy never gets stretched, or stretches slowly.  A reducible or
finite-order braid would indicate poor design.  The entropy is a measure of
the taffy puller's effectiveness: it gives the rate of growth of curves
anchored on the rods.  Thus, the length of the taffy is multiplied
(asymptotically) by $\ee^{1.7627} \simeq 5.828$ for each full period of rod
motion.  Needless to say, this leads to extremely rapid growth, since after
$10$ periods the taffy length has been multiplied by roughly $10^7$.

The design in Figure~\ref{fig:taffy_4rods} can be plotted and analyzed with
\begin{lstlisting}[frame=single,framerule=0pt]
>> b = taffy('4rods')

b = < 1  3  2  2  1  3 >
\end{lstlisting}
When we apply \lstinline{train} to this braid we find the braid is
pseudo-Anosov with exactly the same entropy as the 3-rod taffy puller,
$1.7627$.  There is thus no obvious advantage to using more rods in this case.

\begin{figure}
\begin{center}
\subfigure[]{
  \includegraphics[height=.2\textheight]{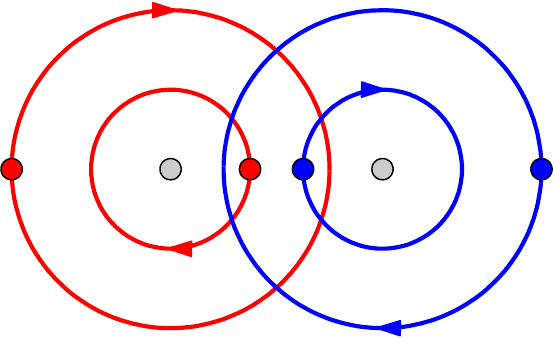}
  \label{fig:taffy_6rods-bad}
}\hspace{1em}
\subfigure[]{
  \includegraphics[height=.2\textheight]{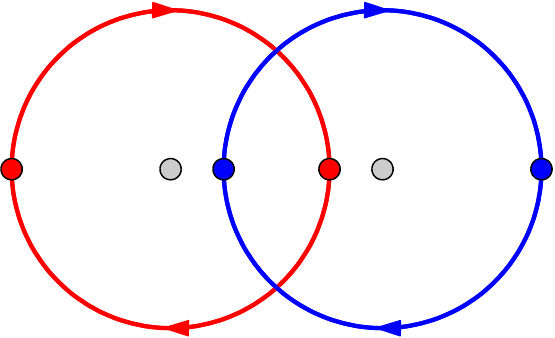}
  \label{fig:taffy_6rods}
}
\end{center}
\caption{(a) A six-rod taffy puller based on Figure~\ref{fig:taffy_4rods},
  with two added fixed rods (gray).  This is a poor design, since it leads to
  a reducible braid. %
  \index{braid!reducible}%
  (b) Same as (a), but with the same radius of motion for all the moving rods.
  The braid is in this case pseudo-Anosov, with larger entropy than the 4-rod
  design.}
\label{fig:taffy_6rods-improved}
\end{figure}

A simple modification of the 4-rod design in Figure~\ref{fig:taffy_4rods} is
shown in Figure~\ref{fig:taffy_6rods-bad}.
The only change is to extend the rotation axles into two extra fixed rods
(shown in gray).  The resulting braid is
\begin{lstlisting}[frame=single,framerule=0pt]
>> b = taffy('6rods-bad')

b = < 2  1  2  4  5  4  3  3  2  1  2  4  5  4 >
\end{lstlisting}
with Thurston--Nielsen type
\index{braid!reducible|(}%
\begin{lstlisting}[frame=single,framerule=0pt]
>> train(b)

ans = struct with fields:

       braid: [1x1 braidlab.braid]
      tntype: 'reducible'
     entropy: 0
    transmat: [5x5 double]
       ttmap: {11x1 cell}
\end{lstlisting}
There are reducing curves in this design: simply wrap a loop around the left
gray rod and the inner red rod, and it will rotate without stretching. %
\index{braid!reducible|)}%
To avoid this, we extend the radius of motion of the inner rods to equal that
of the outer ones, and obtain the design shown in
Figure~\ref{fig:taffy_4rods}.  The corresponding braid is
\begin{lstlisting}[frame=single,framerule=0pt]
>> b = taffy('6rods')

b = < 3  2  1  2  4  5  4  3  3  2  1  2  5  4  5  3 >
\end{lstlisting}
with Thurston--Nielsen type and entropy
\begin{lstlisting}[frame=single,framerule=0pt]
>> [t,entr] = train(b)

ans = struct with fields:

       braid: [1x1 braidlab.braid]
      tntype: 'pseudo-Anosov'
     entropy: 2.6339
    transmat: [5x5 double]
       ttmap: {11x1 cell}
\end{lstlisting}
\index{braid class@\braid\ class!train@\lstinline{train}|)}%
The fixed rods have increased the entropy by~$50\%$!  This sounds like a
fairly small change, but what it means is that this 6-rod design achieves
growth of~$10^7$ in about~$6$ iterations rather than~$10$.  Alexander Flanagan
constructed this six-rod device while an undergraduate student at the
University of Wisconsin -- Madison, but as far as we know this new design has
not yet been used in commercial applications.

The symmetric design of the taffy pullers illustrates one pitfall when
constructing braids.  If we give an optional projection angle %
\index{projection line!bad choice of angle|(}%
of~$\pi/2$ to \lstinline{taffy}:
\begin{lstlisting}[frame=single,framerule=0pt]
>> taffy('4rods',pi/2)
Error using braidlab.braid/colorbraiding
Paths of particles 2 and 1 have a coincident projection.
Try changing the projection angle.
\end{lstlisting}
This corresponds to using the $y$ (vertical) axis to compute the braid, but as
we can see from Figure~\ref{fig:taffy_4rods} this is a bad choice, since all
the rods are initially perfectly aligned along that axis.  The braid obtained
would depend sensitively on numerical roundoff when comparing the rod
projections.%
\footnote{\braidlab\ uses a property \lstinline{BraidAbsTol} to determine when
  coordinates are close enough to be considered coincident, with a default
  value of~$10^{-10}$.  See Section~\ref{sec:prop} for how to set global
  properties.}
Instead of attempting to construct the braid, \braidlab\ returns an error and
asks the user to modify the projection axis.  A tiny change in the projection
line is sufficient to break the symmetry:
\begin{lstlisting}[frame=single,framerule=0pt]
>> taffy('4rods',pi/2 + .01)

ans = < -2  2  1  3  2 -3 -1  3  1  2  1  3 >

>> compact(ans)

ans = < 3  1  2  2  3  1 >
\end{lstlisting}
\index{taffy@\lstinline{taffy}|)}%
which is actually equal to the braid formed from projecting on the~$x$ axis,
though it need only be conjugate %
\index{braid!conjugate}%
(see Section~\ref{sec:braidfromdata}). %
\index{projection line!bad choice of angle|)}%

\index{taffy pullers|)}

\section{Side note: On filling-in punctures}

\index{punctures!filling-in|(}

Recall the command~\lstinline{subbraid}
\index{braid class@\braid\ class!subbraid class@\lstinline{subbraid}}%
from Section~\ref{sec:braidclass}.  We
took the~$4$-string braid~$\sigma_1\sigma_2\sigma_3^{-1}$ and discarded the
third string, to obtain~$\sigma_1\sigma_2^{-1}$:
\begin{lstlisting}[frame=single,framerule=0pt]
>> a = braid([1 2 -3]);
>> b = subbraid(a,[1 2 4])   % discard string 3, keep 1,2,4

b = < 1 -2 >
\end{lstlisting}
\begin{figure}
\begin{center}
\subfigure[]{
  \includegraphics[width=.22\textwidth]{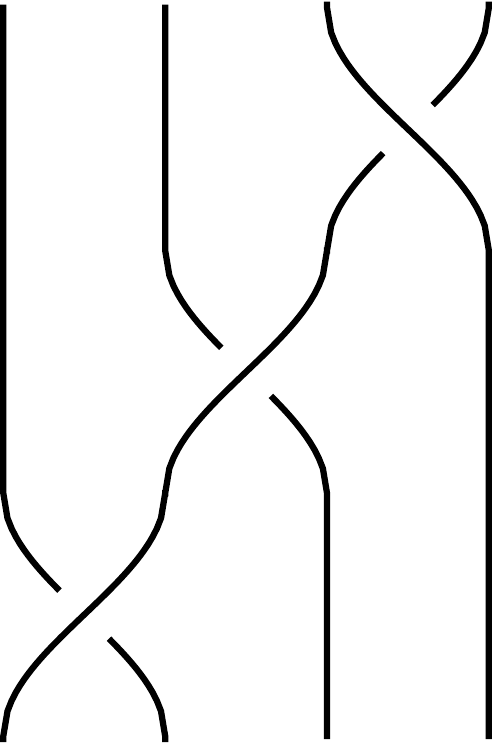}
  \label{fig:s1s2s-3_diagram}
}\hspace{5em}
\subfigure[]{
  \includegraphics[width=.22\textwidth]{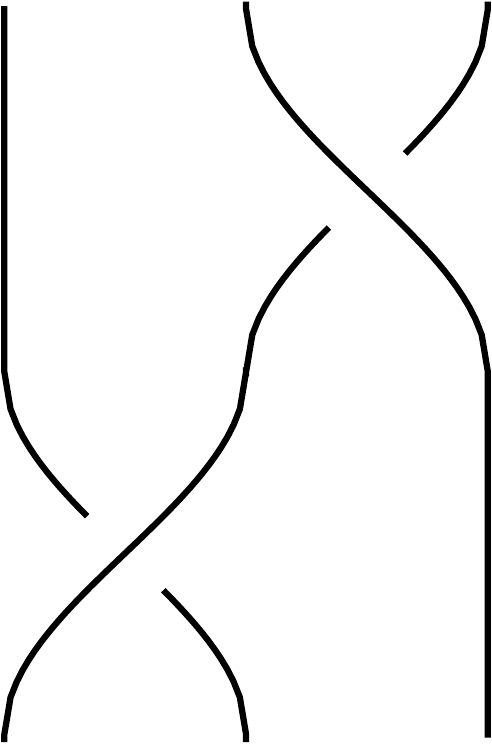}
  \label{fig:s1s-2_diagram}
}
\end{center}
\caption{Removing the third string from the braid
  (a)~$\sigma_1\sigma_2\sigma_3^{-1}$ yields the braid
  (b)~$\sigma_1\sigma_2^{-1}$.}
\label{fig:subbraid}
\end{figure}
The braids \lstinline{a} and \lstinline{b} are shown in
Fig.~\ref{fig:subbraid}; their entropies are %
\index{braid class@\braid\ class!entropy@\lstinline{entropy}|(}%
\index{braid!entropy|(}%
\begin{lstlisting}[frame=single,framerule=0pt]
>> a.entropy, b.entropy

ans = 0.8314
ans = 0.9624
\end{lstlisting}
Note that the entropy of the subbraid~\lstinline{b} is \emph{higher} than the
original braid.  This is counter-intuitive: shouldn't removing strings cause
loops to shorten, therefore lowering their growth?\footnote{In fact, the
  entropy obtained by the removal of a string is constrained by the minimum
  possible entropy %
  \index{braid!entropy!minimum}%
  for the remaining number of strings
  \citep{Song2002,Hironaka2006,Thiffeault2006,
    Ham2007,Venzke_thesis,LanneauThiffeault2011_braids}.  So here the entropy
  of the 3-braid could only be zero or $\ge 0.9624$.}

\begin{figure}
\begin{center}
\subfigure[]{
  \includegraphics[height=.3\textheight]{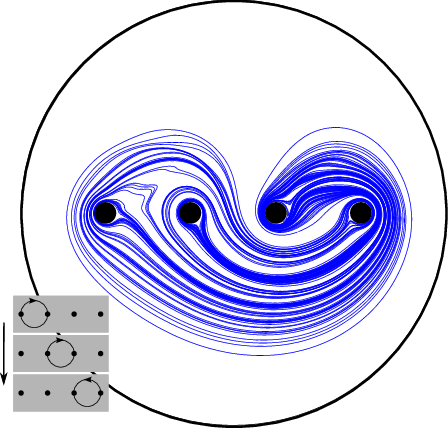}
  \label{fig:s1s2s-3_no_text}
}\hspace{2em}
\subfigure[]{
  \includegraphics[height=.3\textheight]{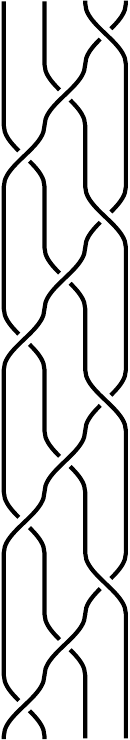}
  \label{fig:s1s2s-3_4_diagram}
}\hspace{2em}
\subfigure[]{
  \includegraphics[height=.3\textheight]{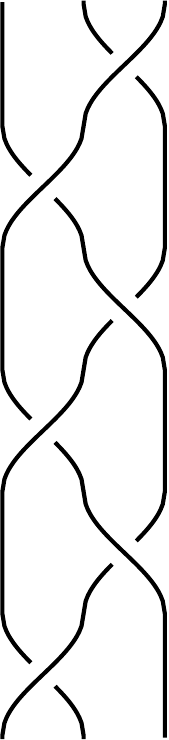}
  \label{fig:s1s-2s1s-2s1s2_diagram}
}
\end{center}
\caption{(a) The mixing protocol specified by the
  braid~$\sigma_1\sigma_2\sigma_3^{-1}$ \citep{Thiffeault2008b}.  The inset
  shows how the rods are moved.  (b)~The pure \index{braid!pure}
  braid~$(\sigma_1\sigma_2\sigma_3^{-1})^4$.  (c)~The
  braid~$(\sigma_1\sigma_2^{-1})^2\sigma_1\sigma_2$, obtained by removing the
  third string from~(b).}
\end{figure}

In some sense this must be true: consider the rod-stirring device shown in
Fig.~\ref{fig:s1s2s-3_no_text}, where the rods move according the to
braid~$\sigma_1\sigma_2\sigma_3^{-1}$.
Removing the third string can be regarded as \emph{filling-in} the third
puncture (rod); clearly then the material line can be shortened, leading to a
decrease in entropy.

\index{braid class@\braid\ class!perm@\lstinline{perm}|(}%
The flaw in the argument is that even though we can remove any string, we
cannot fill in a puncture that is permuted, since the resulting braid does not
define a homeomorphism on the filled-in surface.  To remedy this, let us take
enough powers of the braid~$\sigma_1\sigma_2\sigma_3^{-1}$ to ensure that the
third puncture returns to its original position, using the method
\lstinline{perm} to find the permutation induced by the braid:
\begin{lstlisting}[frame=single,framerule=0pt]
>> perm(a)

ans = 2     3     4     1
\end{lstlisting}
The permutation is cyclic (it can be constructed with exactly one cycle), so
the fourth power should do it:
\begin{lstlisting}[frame=single,framerule=0pt]
>> perm(a^4)

ans = 1     2     3     4
\end{lstlisting}
This is now a pure braid: all the strings return to their original position
(Fig.~\ref{fig:s1s2s-3_4_diagram}).  Now here's the surprise: the subbraid
obtained by removing the third string from \lstinline{a^4} is
\begin{lstlisting}[frame=single,framerule=0pt]
>> b2 = subbraid(a^4,[1 2 4])

b2 = < 1 -2  1 -2  1  2 >
\end{lstlisting}
which is \emph{not} \lstinline{b^4} (Fig.~\ref{fig:s1s-2s1s-2s1s2_diagram})!
However, now there is no paradox in the entropies:\footnote{\citet{Song2005}
  showed that the entropy of a pure braid \index{braid!pure} is greater
  than~$\log(2+\sqrt5) \simeq 1.4436$, if it is nonzero.}
\index{braid!finite-order|(}%
\begin{lstlisting}[frame=single,framerule=0pt]
>> entropy(a^4), entropy(b2)

ans = 3.3258

Warning: Failed to converge to requested tolerance; braid is likely finite-order or has low entropy.  Returning zero entropy.

ans = 0
\end{lstlisting}
\braidlab\ has trouble computing the entropy because the braid \lstinline{b2}
appears to be finite-order.  Indeed, the braid \lstinline{b2} is conjugate
to~$\sigma_1^2$: %
\index{braid class@\braid\ class!compact@\lstinline{compact}}%
\begin{lstlisting}[frame=single,framerule=0pt]
>> c = braid([2 -1],3);
>> compact(c*b2*c^-1)

ans = < 1  1 >
\end{lstlisting}
showing that its entropy is indeed zero.
\index{braid!finite-order|)}%

The moral is: when filling-in punctures, make sure that the strings being
removed are permuted only among themselves.  For very long, random braids, we
still expect that removing a string will decrease the entropy, since the
string being removed will have returned to its initial position many times. %
\index{braid class@\braid\ class!entropy@\lstinline{entropy}|)}%
\index{braid!entropy|)}%
\index{braid class@\braid\ class!perm@\lstinline{perm}|)}%
\index{punctures!filling-in|)}

\section{Setting global properties}
\label{sec:prop}

Braids have been studied for a long time and by many communities, so several
different conventions were bound to emerge.  \braidlab\ has some reasonable
default conventions, but perhaps you'd be happier if it used your favorite
one.  Luckily, \braidlab\ has you covered.

\begin{figure}
  \begin{center}
    \includegraphics[width=\textwidth]{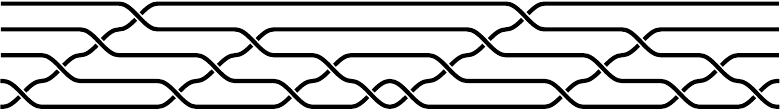}
  \end{center}
  \caption{\lstinline{plot(braid('FullTwist',5))} with the property
    \lstinline{'BraidPlotDir'} set to \hbox{\lstinline{'lr'}}.}
  \label{fig:fulltwist_lr}
\end{figure}

\index{prop@\lstinline{prop}|(}%
To see the properties available and their current values, use the
\lstinline{prop} command:
\begin{lstlisting}[frame=single,framerule=0pt]
>> prop

ans =         GenRotDir: 1
          GenLoopActDir: 'lr'
       GenPlotOverUnder: 1
            BraidAbsTol: 1.0000e-10
           BraidPlotDir: 'bt'
    LoopCoordsBasePoint: 'right'
\end{lstlisting}
To set a property, use something like \lstinline{prop('BraidPlotDir','lr')}.
This will plot braids from left-to-right from now on, as in
Fig.~\ref{fig:fulltwist_lr}.
The default is \lstinline{'bt'} for bottom-to-top.  Other possible values are
\lstinline{'tb'} for top-to-bottom and \lstinline{'rl'} for right-to-left.
See \lstinline{help prop} for more information on properties and allowed
values.
\index{prop@\lstinline{prop}|)}

\section*{Acknowledgments}
\addcontentsline{toc}{section}{Acknowledgments}

The development of \braidlab\ was supported by the US National Science
Foundation, under grants DMS-0806821 and CMMI-1233935. The authors thank
Michael Allshouse and Margaux Filippi for extensive testing and comments.
Michael Allshouse also contributed some of the code. Alexander Flanagan helped
with testing and the taffy puller research. James Puckett and Karen Daniels
provided the test data from their granular medium
experiments~\citep{Puckett2012}; \braidlab\ uses Toby Hall's \emph{Train}
\citep{HallTrain}; Jae Choon Cha's \emph{CBraid} \citep{CBraid}; Juan
Gonz\'{a}lez-Meneses's \emph{Braiding} \citep{Braiding}; John D'Errico's
\emph{Variable Precision Integer Arithmetic} \citep{vpi}; Markus Buehren's
\emph{assignmentoptimal} \citep{assignmentoptimal}; and Jakob Progsch's
\emph{ThreadPool} \citep{ThreadPool}.

\appendix

\section{Installing \braidlab}
\label{sec:install}
\index{installing \braidlab}

\index{Matlab!MEX files|(} \braidlab\ consists of Matlab files together with C
and C++ auxiliary files, so-called MEX files.  The MEX files are used to
greatly speed up calculations.  Many commands will work even if the MEX files
are unavailable, but much more slowly.  (A few commands won't work at all.)
However, MEX files need to be first compiled with Matlab's \lstinline{mex}
compiler.  \index{Matlab!MEX files|)}

\subsection{Precompiled packages}
\label{sec:precompiled}

Some tar files of the precompiled latest released version are
available at \url{http://github.com/jeanluct/braidlab/releases}.  If one of
those suits your system, then download and untar:
\begin{lstlisting}[frame=single,framerule=0pt,escapechar=*,%
  language=bash,backgroundcolor=\color{white}]
$ tar xvzf braidlab-<version>.tar.gz
\end{lstlisting}
A binary might still work even if the operating system and Matlab versions
don't match exactly.

For Mac OS~X, there is a version compiled without
\href{https://gmplib.org/}{GMP} (the GNU MultiPrecision library) %
\index{GMP}%
since it is a bit tricky to install (see Section~\ref{sec:gmp}).

On Ubuntu or Debian Linux, you can install the GMP source and binaries with
\begin{lstlisting}[frame=single,framerule=0pt,escapechar=*,%
  language=bash,backgroundcolor=\color{white}]
$ sudo apt-get install libgmp-dev
\end{lstlisting}
After you install the binaries you must set the Matlab path as described in
Section~\ref{sec:path}.

\index{OS X|(}%
Mac OS~X does not have access to \lstinline{apt} or
\lstinline{aptitude}. Their functionality is in part reproduced by the %
\index{Homebrew|see{OS X}}%
\index{OS X!Homebrew}%
\href{http://brew.sh/}{Homebrew project}. If you install Homebrew, you can
easily obtain GMP by running
\begin{lstlisting}[frame=single,framerule=0pt,escapechar=*,%
  language=bash,backgroundcolor=\color{white}]
$ brew install gmp
\end{lstlisting}
Since \lstinline{brew} is conservative about installing libraries and
sources into system directories, you will likely need to set
\lstinline{LD_LIBRARY_PATH}, and \lstinline{PATH} variables to include
the location of the local Homebrew library storage (see Homebrew
documentation for further instructions).

\index{OS X|)}%

\subsection{Cloning the repository}
\label{sec:cloning}

\index{Git} If you prefer to have the latest (possibly unstable) development
version, and know how to compile Matlab MEX %
\index{Matlab!MEX files|(} files on your system, then you can clone the GitHub
source repository with the terminal command
\begin{lstlisting}[frame=single,framerule=0pt,escapechar=*,%
  language=bash,backgroundcolor=\color{white}]
$ git clone https://github.com/jeanluct/braidlab.git
\end{lstlisting}
assuming Git is installed on your system.  If you prefer to use
Mercurial\index{Mercurial}, make sure you have the
\href{http://hg-git.github.io/}{\lstinline{hg-git}} extension enabled and type
\begin{lstlisting}[frame=single,framerule=0pt,escapechar=*,%
  language=bash,backgroundcolor=\color{white}]
$ hg clone git+https://github.com/jeanluct/braidlab.git
\end{lstlisting}
Either way, after the cloning finishes type
\begin{lstlisting}[frame=single,framerule=0pt,escapechar=*,%
  language=bash,backgroundcolor=\color{white}]
$ cd braidlab; make
\end{lstlisting}
to compile the MEX files.  (If for some reasone you need to recompile, run
\lstinline{make clean} first.)  Note that you can still use \braidlab\ even if
you're unable to compile the MEX files, but some commands will be unavailable
or run (much) more slowly. %
\index{Matlab!MEX files|)}

\index{GMP|(}%
If you receive error messages because \href{https://gmplib.org/}{GMP} (the GNU
MultiPrecision library) is not installed on your system, instead of the above
use
\begin{lstlisting}[frame=single,framerule=0pt,escapechar=*,%
  language=bash,backgroundcolor=\color{white}]
$ cd braidlab; make BRAIDLAB_USE_GMP=0
\end{lstlisting}
This will slow down some functions, in particular testing for equality of
large braids.  If you want to install GMP first, see Section~\ref{sec:gmp}.

\index{GMP|)}

\subsection{Setting Matlab's path}
\label{sec:path}

\lstset{language=Matlab}
\lstset{breaklines=true}
\lstset{backgroundcolor=\color{beige}}

The package \braidlab\ is defined inside a Matlab namespace,
\index{Matlab!namespace} which are specified as subfolders beginning with a
`\lstinline{+}' character.  The Matlab path \index{Matlab!path} must contain
the folder that contains the subfolder \lstinline{+braidlab}, and not the
\lstinline{+braidlab} folder itself:
\begin{lstlisting}[frame=single,framerule=0pt,escapechar=*]
>> addpath '*\it path to folder containing +braidlab*'
\end{lstlisting}
To execute a \braidlab\ \textit{function}, either call it using the syntax
\hbox{\lstinline{braidlab.}\textit{function}}, or import the whole namespace:
\begin{lstlisting}[frame=single,framerule=0pt]
>> import braidlab.*
\end{lstlisting}
This allows invoking \textit{function} by itself, without the
\lstinline{braidlab} prefix.  For the remainder of this document, we
assume this has been done and omit the \lstinline{braidlab} prefix.
The \lstinline{addpath} and \lstinline{import} commands can be added
to \lstinline{startup.m} to ensure they are executed at the start of
every Matlab session.

\subsection{Testing your installation}

\index{testsuite}
To check that everything is working, \braidlab\ includes a testsuite.  From
Matlab, change to the \lstinline{testsuite} folder, and run
\begin{lstlisting}[frame=single,framerule=0pt]
>> test_braidlab
\end{lstlisting}
making sure the path is set properly (Section~\ref{sec:path}).
Note that running the testsuite requires Matlab version 2013a or later.

\subsection{Troubleshooting the installation process}
\label{sec:trouble}
\index{installing \braidlab!troubleshooting|(}
\index{troubleshooting!installation|(}

Here are some common problems that can occur when installing or compiling
braidlab.

\subsubsection{Unsupported compiler}

\index{Linux|(}%
Linux distributions often use very recent C/C++ compilers that are not yet
supported by Matlab.  If you get such an error from MEX, %
\index{Matlab!MEX files}%
it will tell you which version of %
\index{GCC compiler|(}%
GCC it wants.  For example, if it claims it needs GCC 4.7 or earlier, you can
try
\begin{lstlisting}[frame=single,framerule=0pt,escapechar=*,%
  language=bash,backgroundcolor=\color{white}]
$ which gcc-4.7
\end{lstlisting}
to see if a path to the command exists.  If it does, you have an earlier
compiler installed and you can skip to the next paragraph.  If the
\lstinline{which} command above didn't return anything, you can try to install
an older version of GCC:
\begin{lstlisting}[frame=single,framerule=0pt,escapechar=*,%
  language=bash,backgroundcolor=\color{white}]
$ sudo apt-get install gcc-4.7 g++-4.7
\end{lstlisting}
This last line is for Ubuntu and Debian Linux distributions.  Note that this
will \emph{not} overwrite the default compiler.  If your Linux distribution
doesn't allow you to easily install the required compiler, you could always
\href{https://gcc.gnu.org/releases.html}{compile and install it from scratch!}
That's fairly tedious, though.

Assuming that you now have a Matlab-supported compiler installed, you need to
tell MEX %
\index{Matlab!MEX files}%
to use it.  If you're lucky, this could be as easy as typing
\begin{lstlisting}[frame=single,framerule=0pt,escapechar=*,%
  language=bash,backgroundcolor=\color{white}]
$ mex -setup cpp
\end{lstlisting}
but typically this will just lead to Matlab using the wrong compiler anyways.
In that case you can try and download
\href{http://github.com/jeanluct/braidlab/raw/develop/devel/gcc-alternatives}{this
  script}, then run it:
\begin{lstlisting}[frame=single,framerule=0pt,escapechar=*,%
  language=bash,backgroundcolor=\color{white}]
$ bash ./gcc-alternatives
\end{lstlisting}
You will be prompted for your password (assuming you have administrator
privileges), and then asked to select a compiler version.  You should then be
able to compile \braidlab.  You can later switch back to your original
compiler with the command
\begin{lstlisting}[frame=single,framerule=0pt,escapechar=*,%
  language=bash,backgroundcolor=\color{white}]
$ sudo update-alternatives --config gcc
\end{lstlisting}

\index{GCC compiler|)}%
\index{Linux|)}%

\subsubsection{Compiling with GMP}
\label{sec:gmp}

\index{GMP|(}%
\href{https://gmplib.org/}{GMP} is the GNU MultiPrecision library.  \braidlab\
uses GMP to handle arbitrarily large integers, which is necessary for testing
equality of very long braids.  If GMP is not available, \braidlab\ falls back
on the VPI library \citep{vpi}, which is written natively in Matlab and hence
is much slower.

If you wish to compile \braidlab\ from source using GMP, you first need to
install GMP (see \ref{sec:precompiled}), and then make the GMP headers
\lstinline{gmp.h} and \lstinline{gmpxx.h} available to the C++ compiler that
Matlab is using. Both the %
\index{GCC compiler}%
GCC compiler (more common on Linux %
\index{Linux}%
systems) and %
\index{Clang compiler|see{OS X}}%
\index{OS X!Clang compiler}%
Clang compiler (more common on Mac OS~X systems) %
\index{OS X}%
use the environment variable \lstinline{CPLUS_INCLUDE_PATH} to specify path to
include files. For example, if both header files reside in
\lstinline{/usr/local/include} on your Linux or Mac OS~X system, you would
issue (in \lstinline{bash} shell):
\begin{lstlisting}[frame=single,framerule=0pt,escapechar=*,%
  language=bash,backgroundcolor=\color{white}]
$ export CPLUS_INCLUDE_PATH=$CPLUS_INCLUDE_PATH:/usr/local/include
\end{lstlisting}
or in \lstinline{tcsh} shell (now the default on Mac OS~X):
\begin{lstlisting}[frame=single,framerule=0pt,escapechar=*,%
  language=csh,backgroundcolor=\color{white}]
$ setenv CPLUS_INCLUDE_PATH $CPLUS_INCLUDE_PATH:/usr/local/include
\end{lstlisting}
before issuing the \lstinline{make} command to build \lstinline{braidlab}.

\index{GMP|)}%

\subsubsection{Polish \LaTeX\ gets in the way}

\index{Linux|(}%
This is a strange one.  If on compilation you see an error like this:
\begin{lstlisting}[frame=single,framerule=0pt,escapechar=*,%
  language=bash,backgroundcolor=\color{white}]
mex: unrecognized option `-largeArrayDims'
mex: unrecognized option `-O'
mex: unrecognized option `-DBRAIDLAB_USE_GMP'
This is pdfTeX, Version 3.1415926-2.5-1.40.14 (TeX Live 2013)
 restricted \write18 enabled.
entering extended mode
! I can't find file `"CFLAGS=-O -DMATLAB_MEX_FILE"'.
\end{lstlisting}
This is due to the command \lstinline{mex} --- part of the Polish \LaTeX\
package --- shadowing Matlab's \lstinline{mex} compiler.  A simple solution,
if you don't use the Polish language often, is to simply remove the package:
\begin{lstlisting}[frame=single,framerule=0pt,escapechar=*,%
  language=bash,backgroundcolor=\color{white}]
$ sudo apt-get remove texlive-lang-polish
\end{lstlisting}
This last line is for Ubuntu and Debian Linux distributions.  You can also
manually rename the Polish \lstinline{mex} command to something like
\lstinline{mex.polish}, and then make sure Matlab's \lstinline{mex} is in your
path.

Another solution is to make sure that the Matlab executable directory appears
early in bash's path variable. For Matlab R2014b on Mac OS~X %
\index{OS X}%
and under \lstinline{bash} shell this reads
\begin{lstlisting}[frame=single,framerule=0pt,escapechar=*,%
  language=bash,backgroundcolor=\color{white}]
$ export PATH=/Applications/MATLAB_R2014b.app/bin:$PATH
\end{lstlisting}
or under \lstinline{tcsh} shell (now the default in Mac OS~X)
\begin{lstlisting}[frame=single,framerule=0pt,escapechar=*,%
  language=bash,backgroundcolor=\color{white}]
$ setenv PATH /Applications/MATLAB_R2014b.app/bin:$PATH
\end{lstlisting}

\index{Linux|)}%

\subsubsection{\lstinline{largeArraydims} warning}

You might get this warning:
\begin{lstlisting}[frame=single,framerule=0pt,escapechar=*,%
  language=bash,backgroundcolor=\color{white}]
Warning: Legacy MEX infrastructure is provided for compatibility; it will be removed in a future version of MATLAB.
\end{lstlisting}
This can be safely ignored.  Matlab is transitioning from a shorter to a
longer type of internal array indexing.  Eventually the
\lstinline{-largeArraydims} flag will be removed from \braidlab.

\index{installing \braidlab!troubleshooting|)}
\index{troubleshooting!installation|)}

\section{Troubleshooting \braidlab}
\label{sec:troubleshooting}

\index{troubleshooting|(}%
If \braidlab\ is behaving unexpectedly, or reporting an error that you cannot
diagnose, tips in this section will help you identify the cause of the problem
and alert the developers if an issue exists.

\subsection{Global flags}
\label{sec:global-flags}

\index{global flags|(}%
There are several flags that you can use to help the identify the root of the
problem:
\begin{itemize}
\item \lstinline{BRAIDLAB_braid_nomex}%
  \index{global flags!BRAIDLABbraidnomex@\lstinline{BRAIDLAB_braid_nomex}}:
  set to \lstinline{true} to \emph{disable} MEX %
  \index{Matlab!MEX files}%
  (C++) versions of algorithms that convert trajectories to a
  \lstinline{braid} or \lstinline{databraid}.
\item \lstinline{BRAIDLAB_loop_nomex}%
  \index{global flags!BRAIDLABloopnomex@\lstinline{BRAIDLAB_loop_nomex}}:
  set to \lstinline{true} to \emph{disable} MEX %
  \index{Matlab!MEX files}%
  (C++) versions of algorithms that apply braid generators to loops.
\item \lstinline{BRAIDLAB_threads}%
  \index{global flags!BRAIDLABthreads@\lstinline{BRAIDLAB_threads}}: set to
  \lstinline{1} to disable parallel MEX %
  \index{Matlab!MEX files}%
  algorithms, or to a larger integer to fix the number of parallel processing
  threads used.
\item \lstinline{BRAIDLAB_debuglvl}%
  \index{global flags!BRAIDLABdebuglvl@\lstinline{BRAIDLAB_debuglvl}}: set
  to \lstinline{1}, \lstinline{2}, \lstinline{3}, \lstinline{4} to
  progressively increase the amount of output produced by
  \braidlab. \emph{Warning:} setting this flag to any positive number may
  result in an overwhelming amount of output lines.
\end{itemize}

All the flags are set through MATLAB's \lstinline{global} variables as follows.
To set the global value of the flag \lstinline{BRAIDLAB_<flagname>} to \lstinline{true}, issue the following set of commands \emph{before} you run your code.
\begin{lstlisting}[frame=single,framerule=0pt]
>> global BRAIDLAB_<flagname>
>> BRAIDLAB_<flagname> = true;
\end{lstlisting}
At this point, any change in the flag value, e.g.,
\lstinline{BRAIDLAB_<flagname> = false}, will be reflected in \braidlab's
operation. When you are done troubleshooting the code, clear the flag by
calling
\begin{lstlisting}[frame=single,framerule=0pt]
>> clear global BRAIDLAB_<flagname>
\end{lstlisting}
to restore the default state of \braidlab.
\index{global flags|)}%

\subsection{Reporting issues and suggestions}
\label{sec:reporting-issues}

We appreciate your feedback on how to improve \braidlab. Before reporting a
bug, please make sure that you are using the most recent version, as the bugs
are continually fixed and new features added.

To inform the developers of a problem with \braidlab, unexpected behavior, or
a suggestion for an improvement, please use the interface on the GitHub
repository:
\begin{center}\url{http://github.com/jeanluct/braidlab/issues} \end{center}

There you can view the list of currently known issues and find out if anyone
else has observed the same problem as you. If not, submit the new issue and
include:
\begin{itemize}
\item The version of \braidlab. If you are updating \braidlab{} using
  \lstinline{git}, please copy-paste the output of \lstinline{git log -1} run
  in your \braidlab{} folder. If you have downloaded a compiled version of
  \braidlab{}, please let us know the version number on this guide's front
  page.
\item A minimal example re-creating the problem. You can also attach any data
  files to the GitHub issue that would help us recreate the problem.
\item Output of the Matlab command window during the behavior of the
  error. Two Matlab features are particularly useful for this purpose: the
  \lstinline{diary} command, which saves the output of all commands between
  \lstinline{diary <filename>} and \lstinline{diary off} to a text file, and
  \lstinline{[msg,id] = lasterr} and \lstinline{[msg,id] = lastwarn}, which
  return the text and the identifier of the last error/warning issued by
  Matlab.
\item Any other information about what you are trying to achieve: before we
  fix the problem, we may be able to suggest a way to circumvent the bug.
\end{itemize}

We will also respond to your requests by e-mail; however, using the GitHub
issues interface helps us keep track of issues and their resolution more
reliably.  \index{troubleshooting|)}

\bibliographystyle{jfm}
{\small
\bibliography{braidlab_guide}

\begin{thebibliography}{42}
\expandafter\ifx\csname natexlab\endcsname\relax\def\natexlab#1{#1}\fi

\bibitem[Artin(1947)]{Artin1947}
{\sc Artin, E.} 1947 Theory of braids. {\em Ann. Math.\/} {\bf 48}~(1),
  101--126.

\bibitem[Bangert {\em et~al.\/}(2002)Bangert, Berger \& Prandi]{Bangert2002}
{\sc Bangert, P.~D., Berger, M.~A. \& Prandi, R.} 2002 In search of minimal
  random braid configurations. {\em J. Phys. A\/} {\bf 35}~(1), 43--59.

\bibitem[Bestvina \& Handel(1995)]{Bestvina1995}
{\sc Bestvina, M. \& Handel, M.} 1995 Train-tracks for surface homeomorphisms.
  {\em Topology\/} {\bf 34}~(1), 109--140.

\bibitem[Bigelow(2001)]{Bigelow2001}
{\sc Bigelow, S.~J.} 2001 Braid groups are linear. {\em J. Amer. Math. Soc.\/}
  {\bf 14}~(2), 471--486.

\bibitem[Birman(1975)]{Birman1975}
{\sc Birman, J.~S.} 1975 {\em Braids, Links, and Mapping Class Groups\/}. {\em
  Annals of Mathematics Studies\/} 82. Princeton, NJ: Princeton University
  Press.

\bibitem[Birman \& Brendle(2005)]{Birman2005}
{\sc Birman, J.~S. \& Brendle, T.~E.} 2005 Braids: {A} survey. In {\em Handbook
  of Knot Theory\/} (ed. W.~Menasco \& M.~Thistlethwaite), pp. 19--104.
  Amsterdam: Elsevier, available at {\tt http://arXiv.org/abs/math.GT/0409205}.

\bibitem[Boyland(1994)]{Boyland1994}
{\sc Boyland, P.~L.} 1994 Topological methods in surface dynamics. {\em
  Topology Appl.\/} {\bf 58}, 223--298.

\bibitem[Budi{\v s}i{\'c} \& {Thiffeault}(2015)]{Budisic2015}
{\sc Budi{\v s}i{\'c}, M. \& {Thiffeault}, J.-L.} 2015 Finite-time braiding
  exponents. {\em Chaos: {An} {Interdisciplinary} {Journal} of {Nonlinear}
  {Science}\/} {\bf 25}~(8), \url{http://arxiv.org/abs/1502.02162}.

\bibitem[Buerhen(2011)]{assignmentoptimal}
{\sc Buerhen, M.} 2011 Functions for the rectangular assignment problem.
  \url{http://www.mathworks.com/matlabcentral/fileexchange/6543}.

\bibitem[Burau(1936)]{Burau1936}
{\sc Burau, W.} 1936 \"{U}ber {Z}opfgruppen und gleichsinnig verdrilte
  {V}erkettungen. {\em Abh. Math. Semin. Hamburg Univ.\/} {\bf 11}, 171--178.

\bibitem[Casson \& Bleiler(1988)]{Casson1988}
{\sc Casson, A.~J. \& Bleiler, S.~A.} 1988 {\em Automorphisms of surfaces after
  {N}ielsen and {T}hurston\/}, {\em London Mathematical Society Student
  Texts\/}, vol.~9. Cambridge: Cambridge University Press.

\bibitem[Cha(2011)]{CBraid}
{\sc Cha, J.~C.} 2011 \textit{CBraid: {A} {C++} library for computations in
  braid groups}. \url{https://github.com/jeanluct/cbraid}.

\bibitem[Dehornoy(2008)]{Dehornoy2008}
{\sc Dehornoy, P.} 2008 Efficient solutions to the braid isotopy problem. {\em
  Discr. Applied Math.\/} {\bf 156}, 3091--3112.

\bibitem[D'Errico(2013)]{vpi}
{\sc D'Errico, J.} 2013 {Variable Precision Integer Arithmetic}.
  \url{http://www.mathworks.com/matlabcentral/fileexchange/22725-variable-precision-integer-arithmetic}.

\bibitem[Dynnikov(2002)]{Dynnikov2002}
{\sc Dynnikov, I.~A.} 2002 On a {Yang--Baxter} map and the {D}ehornoy ordering.
  {\em Russian Math. Surveys\/} {\bf 57}~(3), 592--594.

\bibitem[Dynnikov \& Wiest(2007)]{Dynnikov2007}
{\sc Dynnikov, I.~A. \& Wiest, B.} 2007 On the complexity of braids. {\em
  Journal of the European Mathematical Society\/} {\bf 9}~(4), 801--840.

\bibitem[Fathi {\em et~al.\/}(1979)Fathi, Laundenbach \&
  Po\'{e}naru]{Fathi1979}
{\sc Fathi, A., Laundenbach, F. \& Po\'{e}naru, V.} 1979 Travaux de {T}hurston
  sur les surfaces. {\em Ast\'{e}risque\/} {\bf 66-67}, 1--284.

\bibitem[Finn \& Thiffeault(2011)]{MattFinn2011_silver}
{\sc Finn, M.~D. \& Thiffeault, J.-L.} 2011 Topological optimization of
  rod-stirring devices. {\em SIAM Rev.\/} {\bf 53}~(4), 723--743.

\bibitem[Finn {\em et~al.\/}(2006)Finn, Thiffeault \& Gouillart]{MattFinn2006}
{\sc Finn, M.~D., Thiffeault, J.-L. \& Gouillart, E.} 2006 Topological chaos in
  spatially periodic mixers. {\em Physica D\/} {\bf 221}~(1), 92--100.

\bibitem[Gonz\'{a}lez-Meneses(2011)]{Braiding}
{\sc Gonz\'{a}lez-Meneses, J.} 2011 \textit{Braiding: {A} computer program for
  handling braids}. The version used is distributed with \textit{CBraid}:
  \url{http://code.google.com/p/cbraid}.

\bibitem[Hall(2012)]{HallTrain}
{\sc Hall, T.} 2012 \textit{Train: {A} {C++} program for computing train tracks
  of surface homeomorphisms}. \url{http://www.liv.ac.uk/~tobyhall/T_Hall.html}.

\bibitem[Hall \& Yurtta\c{s}(2009)]{Hall2009}
{\sc Hall, T. \& Yurtta\c{s}, S.~{\"{O}}.} 2009 On the topological entropy of
  families of braids. {\em Topology Appl.\/} {\bf 156}~(8), 1554--1564.

\bibitem[Ham \& Song(2007)]{Ham2007}
{\sc Ham, J.-Y. \& Song, W.~T.} 2007 The minimum dilatation of pseudo-{A}nosov
  5-braids. {\em Experiment. Math.\/} {\bf 16}~(2), 167--179.

\bibitem[Hironaka \& Kin(2006)]{Hironaka2006}
{\sc Hironaka, E. \& Kin, E.} 2006 A family of pseudo-{A}nosov braids with
  small dilatation. {\em Algebraic \& Geometric Topology\/} {\bf 6}, 699--738.

\bibitem[Kassel \& Turaev(2008)]{KasselTuraev}
{\sc Kassel, C. \& Turaev, V.} 2008 {\em Braid groups\/}. New York, NY:
  Springer.

\bibitem[Lanneau \& Thiffeault(2011)]{LanneauThiffeault2011_braids}
{\sc Lanneau, E. \& Thiffeault, J.-L.} 2011 On the minimum dilatation of braids
  on the punctured disc. {\em Geometriae Dedicata\/} {\bf 152}~(1), 165--182.

\bibitem[Lawrence(1990)]{Lawrence1990}
{\sc Lawrence, R.~J.} 1990 Homological representations of the {H}ecke algebra.
  {\em Comm. Math. Phys.\/} {\bf 135}~(1), 141--191.

\bibitem[Moussafir(2006)]{Moussafir2006}
{\sc Moussafir, J.-O.} 2006 On computing the entropy of braids. {\em Func.
  Anal. and Other Math.\/} {\bf 1}~(1), 37--46.

\bibitem[Paterson \& Razborov(1991)]{Paterson1991}
{\sc Paterson, M.~S. \& Razborov, A.~A.} 1991 The set of minimal braids is
  co-{NP} complete. {\em J. Algorithm\/} {\bf 12}, 393--408.

\bibitem[Progsch(2012)]{ThreadPool}
{\sc Progsch, J.} 2012 {ThreadPool}: A simple {C++11} thread pool
  implementation. \url{https://github.com/progschj/ThreadPool}.

\bibitem[Puckett {\em et~al.\/}(2012)Puckett, Lechenault, Daniels \&
  Thiffeault]{Puckett2012}
{\sc Puckett, J.~G., Lechenault, F., Daniels, K.~E. \& Thiffeault, J.-L.} 2012
  Trajectory entanglement in dense granular materials. {\em Journal of
  Statistical Mechanics: Theory and Experiment\/} {\bf 2012}~(6), P06008.

\bibitem[Song(2005)]{Song2005}
{\sc Song, W.~T.} 2005 Upper and lower bounds for the minimal positive entropy
  of pure braids. {\em Bull. London Math. Soc.\/} {\bf 37}~(2), 224--229.

\bibitem[Song {\em et~al.\/}(2002)Song, Ko \& Los]{Song2002}
{\sc Song, W.~T., Ko, K.~H. \& Los, J.~E.} 2002 Entropies of braids. {\em J.
  Knot Th. Ramifications\/} {\bf 11}~(4), 647--666.

\bibitem[Thiffeault(2005)]{Thiffeault2005}
{\sc Thiffeault, J.-L.} 2005 Measuring topological chaos. {\em Phys. Rev.
  Lett.\/} {\bf 94}~(8), 084502.

\bibitem[Thiffeault(2010)]{Thiffeault2010}
{\sc Thiffeault, J.-L.} 2010 Braids of entangled particle trajectories. {\em
  Chaos\/} {\bf 20}, 017516.

\bibitem[Thiffeault(2018)]{Thiffeault2018}
{\sc Thiffeault, J.-L.} 2018 The mathematics of taffy pullers. {\em Math.
  Intelligencer\/} {\bf 40}~(1), 26--35,
  \href{https://arxiv.org/abs/1608.00152}{arXiv:1608.00152}.

\bibitem[Thiffeault \& Finn(2006)]{Thiffeault2006}
{\sc Thiffeault, J.-L. \& Finn, M.~D.} 2006 Topology, braids, and mixing in
  fluids. {\em Phil. Trans. R. Soc. Lond. A\/} {\bf 364}, 3251--3266.

\bibitem[Thiffeault {\em et~al.\/}(2008)Thiffeault, Finn, Gouillart \&
  Hall]{Thiffeault2008b}
{\sc Thiffeault, J.-L., Finn, M.~D., Gouillart, E. \& Hall, T.} 2008 Topology
  of chaotic mixing patterns. {\em Chaos\/} {\bf 18}, 033123.

\bibitem[Thurston(1988)]{Thurston1988}
{\sc Thurston, W.~P.} 1988 On the geometry and dynamics of diffeomorphisms of
  surfaces. {\em Bull. Am. Math. Soc.\/} {\bf 19}, 417--431.

\bibitem[Venzke(2008)]{Venzke_thesis}
{\sc Venzke, R.~W.} 2008 Braid forcing, hyperbolic geometry, and
  pseudo-{A}nosov sequences of low entropy. PhD thesis, California Institute of
  Technology.

\bibitem[Weisstein(2013)]{AlexanderPolynomial}
{\sc Weisstein, E.~W.} 2013 Alexander polynomial. From \textit{MathWorld}---A
  Wolfram Web Resource.
  \url{http://mathworld.wolfram.com/AlexanderPolynomial.html}.

\bibitem[Yurtta\c{s}(2014)]{Yurttas2014_preprint}
{\sc Yurtta\c{s}, S.~{\"{O}}.} 2014 Dynnikov and train track transition
  matrices of pseudo-{A}nosov braids. \url{http://arxiv.org/abs/1402.4378}.

\end{thebibliography}
}
\addcontentsline{toc}{section}{References}


\index{crossing|seealso{projection line}}%
\index{projection line|seealso{crossing}}%
\index{Finite Time Braiding Exponent (FTBE)|seealso{\lstinline{databraid} class}}%

\clearpage
\addcontentsline{toc}{section}{Index}
\printindex

\end{document}